\documentclass{article}
\usepackage{mathrsfs}
\usepackage{latexsym}
\usepackage{graphicx}
\usepackage{amssymb}
\usepackage{amsmath}
\usepackage{calrsfs}

\begin{document}
\title{\normalsize {ON THE REGULARITY FOR
3D NAVIER-STOKES EQUATION } }

\author {\normalsize{Qun Lin}  \footnotemark \thanks{School of Mathematical Sciences, Xiamen University, P.
R. China ( e-mail: \qquad \qquad \qquad linqun@xmu.edu.cn ) } }

\date{\begin{footnotesize}July 12, 2011 \end{footnotesize}}
\maketitle
\bigskip

\begin{minipage}{11cm}\begin{footnotesize} {\bf Abstract.} In this
paper we will prove that the vorticity belongs to
$L^{\infty}(0,T;L^2(\Omega))$ for 3D incompressible Navier-Stokes
equation with periodic initial-boundary value conditions, then the
existence of a global smooth solution is obtained. Our approach is
to construct a set of auxiliary problems to approximate the original
one of vorticity equation.


\end{footnotesize}

\end{minipage}

\bigskip
\bigskip
\bigskip

\textbf{1. Introduction}

Let $\Omega =(0,1)^3$, and $\mathscr{D}(\Omega)$ be the space of
$C^{\infty}$ functions with compact support contained in $\Omega$.
Some basic spaces will be used in this paper:
\begin{equation*}
\begin{split}
 &\mathscr{V}=\{u \in \mathscr{D}(\Omega),\;\,\mbox{div}\,u =0\}  \\
 &V=\mbox{the}\;\mbox{closure}\;\mbox{of}\; \mathscr{V}\; \mbox{in} \;
H^1(\Omega )  \\
 &H=\mbox{the}\;\mbox{closure}\;\mbox{of}\; \mathscr{V}\; \mbox{in}\;
L^2(\Omega )  \\
 \end{split}
\end{equation*}

The velocity-pressure form for Navier- Stokes equation is
\begin{equation}
\label{eq1}
\begin{split}
 &\partial _t u_1 +u_1 \partial _{x_1 } u_1 +u_2 \partial _{x_2 } u_1 +u_3
\partial _{x_3 } u_1 +\partial _{x_1 } p=\Delta u_1 \\
 &\partial _t u_2 +u_1 \partial _{x_1 } u_2 +u_2 \partial _{x_2 } u_2 +u_3
\partial _{x_3 } u_2 +\partial _{x_2 } p=\Delta u_2 \\
 &\partial _t u_3 +u_1 \partial _{x_1 } u_3 +u_2 \partial _{x_2 } u_3 +u_3
\partial _{x_3 } u_3 +\partial _{x_3 } p=\Delta u_3 \\
 \end{split}
\end{equation}
with the initial conditions $\left. {(u_1 ,u_2 ,u_3 )} \right|_{t=0}
=(u_{10} ,u_{20} ,u_{30} )(x)$, henceforth we always ignore the assumption of sufficient smoothness of the initial conditions.
Moreover, the periodic boundary conditions are
\[
u_i (x + e_j ,t) = u_i (x,t),\quad \quad i,j = 1,2,3
\]
and the incompressible condition is
\[
\partial _{x_1 } u_1 +\partial _{x_2 } u_2 +\partial _{x_3 } u_3 =0
\]
where $x = (x_1 ,x_2 ,x_3 )$ is a point of ${\mathbb R}^3$, and $e_j
$ is $j^{\mbox{th}}$ unit vector in ${\mathbb R}^3$.
$u=(u_1,u_2,u_3)$ is velocity, $p$ is pressure, and $\nu>0$ is
viscosity.
\\

The vorticity-velocity
form for Navier-Stokes equation is
\begin{equation}
\label{eq4}
\begin{split}
 &\partial _t \omega _1 \,+u_1 \partial _{x_1 } \omega _1 +u_2 \partial _{x_2
} \omega _1 +u_3 \partial _{x_3 } \omega _1 -\omega _1 \partial _{x_1 } u_1
-\omega _2 \partial _{x_2 } u_1 -\omega _3 \partial _{x_3 } u_1 =\Delta
\omega _1 \\
 &\partial _t \omega _2 +u_1 \partial _{x_1 } \omega _2 +u_2 \partial _{x_2 }
\omega _2 +u_3 \partial _{x_3 } \omega _2 -\omega _1 \partial _{x_1 } u_2
-\omega _2 \partial _{x_2 } u_2 -\omega _3 \partial _{x_3 } u_2 =\Delta
\omega _2 \\
 &\partial _t \omega _3 +u_1 \partial _{x_1 } \omega _3 +u_2 \partial _{x_2 }
\omega _3 +u_3 \partial _{x_3 } \omega _3 -\omega _1 \partial _{x_1 } u_3
-\omega _2 \partial _{x_2 } u_3 -\omega _3 \partial _{x_3 } u_3 =\Delta
\omega _3 \\
 \end{split}
\end{equation}
with the initial conditions $\left. {(\omega _1 ,\omega _2 ,\omega _3 )}
\right|_{t=0} =(\omega _{10} ,\omega _{20} ,\omega _{30}
)=(\mbox{curl}u_{10} ,\;\mbox{curl}u_{20} ,\;\mbox{curl}u_{30} )$, and the periodic boundary conditions :
\[
\omega_i (x + e_j ,t) = \omega_i (x,t),\quad \quad i,j = 1,2,3
\]
and the incompressible condition :
\begin{equation*}
\begin{split}
 &\partial _{x_1 } \omega _1 +\partial _{x_2 } \omega _2 +\partial _{x_3 }
\omega _3 =0 \\
 &\partial _{x_1 } u_1 +\partial _{x_2 } u_2 +\partial _{x_3 } u_3 =0 \\
 \end{split}
\end{equation*}
\\

We here recall the global $L^2$-estimate from [4] for the Navier-Stokes equation of velocity-pressure form. Since
\begin{equation*}
\begin{split}
 &\int_\Omega {u_i (u_1 \partial _{x_1 } u_i +u_2 \partial _{x_2 } u_i +u_3
\partial _{x_3 } u_i )} =\frac{1}{2}\int_\Omega {(u_1 \partial _{x_1 } u_i^2
+u_2 \partial _{x_2 } u_i^2 +u_3 \partial _{x_3 } u_i^2 )} \\
 &\quad =-\frac{1}{2}\int_\Omega {u_i^2 (\partial _{x_1 } u_1 +\partial _{x_2 } u_2
+\partial _{x_3 } u_3 )} =0 \qquad \qquad  i=1,2,3\\
 \end{split}
\end{equation*}
\[
\int_\Omega {(u_1 \partial _{x_1 } p+u_2 \partial _{x_2 } p+u_3
\partial _{x_3 } p)} =-\int_\Omega {p\,(\partial _{x_1 } u_1
+\partial _{x_2 } u_2 +\partial _{x_3 } u_3 )} =0  \qquad \quad
\]
and
\begin{equation*}
\begin{split}
 &\int_\Omega {u_i \Delta u_i } =\int_\Omega {u_i (\partial _{x_1 }^2
 u_i
+\partial _{x_2 }^2 u_i +\partial _{x_3 }^2 u_i )}  \\
 &\quad =-\int_\Omega
{((\partial _{x_1 } u_i )^2+(\partial _{x_2 } u_i )^2+(\partial
_{x_3 } u_i )^2)} \\
 \end{split}
\end{equation*}
then
\begin{equation*}
\begin{split}
 &\int_\Omega {u_1 \partial _t \,u_1 } +\int_\Omega {u_1 (u_1 \partial _{x_1
} u_1 +u_2 \partial _{x_2 } u_1 +u_3 \partial _{x_3 } u_1 )} +\int_\Omega
{u_1 \partial _{x_1 } p} =\int_\Omega {u_1 \Delta u_1 } \\
 &\int_\Omega {u_2 \partial _t \,u_2 } +\int_\Omega {u_2 (u_1 \partial _{x_1
} u_2 +u_2 \partial _{x_2 } u_2 +u_3 \partial _{x_3 } u_2 )} +\int_\Omega
{u_2 \partial _{x_2 } p} =\int_\Omega {u_2 \Delta u_2 } \\
 &\int_\Omega {u_3 \partial _t \,u_3 } +\int_\Omega {u_3 (u_1 \partial _{x_1
} u_3 +u_2 \partial _{x_2 } u_3 +u_3 \partial _{x_3 } u_3 )} +\int_\Omega
{u_3 \partial _{x_3 } p} =\int_\Omega {u_3 \Delta u_3 } \\
 \end{split}
\end{equation*}
so that
\begin{equation*}
\begin{split}
 &\frac{1}{2}\partial _t \;\int_\Omega {(u_1^2 +u_2^2 +u_3^2 )}
+\;\int_\Omega {((\partial _{x_1 } u_1 )^2+(\partial _{x_2 } u_1
)^2+(\partial _{x_3 } u_1 )^2+} \\
 &+(\partial _{x_1 } u_2 )^2+(\partial _{x_2 } u_2 )^2+(\partial _{x_3 } u_2
)^2+(\partial _{x_1 } u_3 )^2+(\partial _{x_2 } u_3 )^2+(\partial _{x_3 }
u_3 )^2)=0 \\
 \end{split}
\end{equation*}
it follows that
\[
\int_\Omega {(u_1^2 +u_2^2 +u_3^2 )} +2\;\int_0^T {(\,\left\| {\nabla u_1 }
\right\|_{L^2(\Omega )}^2 +} \left\| {\nabla u_2 } \right\|_{L^2(\Omega )}^2
+\left\| {\nabla u_3 } \right\|_{L^2(\Omega )}^2 )=\int_\Omega {(u_{10}^2
+u_{20}^2 +u_{30}^2 )}
\]

Hence we have
\begin{equation}
\label{eq2} \mathop {\sup }\limits_{t\in (0,T)} \;\;\int_\Omega
{(u_1^2 +u_2^2 +u_3^2 )} <+\infty
\end{equation}
\begin{equation}
\label{eq3} \,\int_0^T {(\,\left\| {\nabla u_1 }
\right\|_{L^2(\Omega )}^2 +} \left\| {\nabla u_2 }
\right\|_{L^2(\Omega )}^2 +\left\| {\nabla u_3 }
\right\|_{L^2(\Omega )}^2 )<+\infty
\end{equation}

Above $u$ can be interpreted as the Galerkin approximation of the
solution, but (3) and (4) are also true for the solution of problem
(1).
\\

The rest of sections are arranged as follows : In section 2 and 3, we introduce a set of auxiliary problems and prove the uniform boundedness and the existence of their solutions in $L^{\infty}(0,T;L^2(\Omega))$. Then it is shown that the solutions of the auxiliary problems converge to that of Naiver-Stokes equation with vorticity-velocity form, which also belongs to $L^{\infty}(0,T;L^2(\Omega))$. Final section will present the solution of Navier-Stokes equation with velocity-pressure form belongs to $L^{\infty}(0,T; H^2(\Omega))$.
\\
\\
\\

\textbf{2. Auxiliary problems}

For the 3D regularity, we only need to prove that the vorticity in (2)
belongs to $L^{\infty}(0,T;L^2(\Omega))$.
\\

Given a partition with respect to $t$ as follows :
\[
0=t_0 <t_1 <t_2 <\cdots <t_{k-1} <t_k <\cdots <t_N =T
\]

On each $t\in (t_{k-1} ,\;t_k )$, we introduce an auxiliary problem :
\begin{equation}
\label{eq5}
\begin{split}
 &\partial _t \tilde {\omega }_1 \,+ \overline{\overline u}_1^k \partial _{x_1 } \overline {\omega
}_1^k + \overline{\overline u}_2^k \partial _{x_2 } \overline
{\omega }_1^k + \overline{\overline u}_3^k
\partial _{x_3 } \overline {\omega }_1^k -\overline {\omega }_1^k \partial
_{x_1 } \overline{\overline u}_1^k -\overline {\omega }_2^k \partial
_{x_2 } \overline{\overline u}_1^k -\overline {\omega }_3^k
\partial
_{x_3 } \overline{\overline u}_1^k +\partial _{x_1 } q=\Delta \tilde {\omega }_1 \\
 &\partial _t \tilde {\omega }_2 + \overline{\overline u}_1^k \partial _{x_1 } \overline {\omega }_2^k
+ \overline{\overline u}_2^k \partial _{x_2 } \overline {\omega
}_2^k + \overline{\overline u}_3^k
\partial _{x_3 } \overline {\omega }_2^k -\overline {\omega }_1^k \partial
_{x_1 } \overline{\overline u}_2^k -\overline {\omega }_2^k \partial
_{x_2 } \overline{\overline u}_2^k -\overline {\omega }_3^k
\partial _{x_3 } \overline{\overline u}_2^k
+\partial _{x_2 } q=\Delta \tilde {\omega }_2 \\
 &\partial _t \tilde {\omega }_3 + \overline{\overline u}_1^k \partial _{x_1 } \overline {\omega }_3^k
+ \overline{\overline u}_2^k \partial _{x_2 } \overline {\omega
}_3^k + \overline{\overline u}_3^k
\partial _{x_3 } \overline {\omega }_3^k -\overline {\omega }_1^k \partial
_{x_1 } \overline{\overline u}_3^k -\overline {\omega }_2^k \partial
_{x_2 } \overline{\overline u}_3^k -\overline {\omega }_3^k
\partial _{x_3 } \overline{\overline u}_3^k
+\partial _{x_3 } q=\Delta \tilde {\omega }_3 \\
 \end{split}
\end{equation}
where the initial value is assumed to be $\tilde {\omega }_i
(x,t_{k-1} )=\tilde {\omega }_i^{k-1} (x) $, \;$\tilde {\omega }_i
(x,0 ) = {\omega }_{i0} (x), \,i=1,2,3 $, \;and

\[
\overline {\omega }_i^k (x)=
 \frac{1}{\Delta t_k }\int_{t_{k-1} }^{t_k } {\tilde {\omega }_i
(x,t)dt}
\]
and
\[
\overline {u }_i^k (x)=
 \frac{1}{\Delta t_k }\int_{t_{k-1} }^{t_k } u_i(x,t)dt, \qquad\;\; i=1,2,3
\]

In addition, let $\varepsilon > 0$, we construct a mollifier
$J_\varepsilon \in C_0^\infty (\mathbb R^3)$ such that

\qquad i) $J_\varepsilon (x) \ge 0,\;\;x \in \mathbb R^3$,

\qquad ii) $J_\varepsilon (x) = 0$ if $\;\left| x \right| \ge
\varepsilon $, and

\qquad iii) $\int_{\mathbb R^3} {J_\varepsilon (x)\, dx} = 1$.

\noindent
\\
then a convolution is defined as
\[
\overline{\overline u} _i^k (x) = J_\varepsilon * \overline u _i^k
(x) = \int_{\mathbb R^3} {J_\varepsilon (x - y)\, \overline u _i^k
(y)\,dy}
\]
where we assume that a zero extension of $\overline u _i^k $ be made
outside $\Omega $.
\\

Similarly we can set the periodic boundary conditions :
\[
\tilde \omega_i (x + e_j ,t) = \tilde \omega_i (x,t),\quad \quad i,j = 1,2,3
\]
and the incompressible condition :
\begin{equation*}
\begin{split}
 &\partial _{x_1 } \tilde \omega _1 + \partial _{x_2 } \tilde \omega _2 + \partial _{x_3 }
\tilde \omega _3 = 0 \\
 &\partial _{x_1 } u_1 +\partial _{x_2 } u_2 +\partial _{x_3 } u_3 =0 \\
 \end{split}
\end{equation*}

It is easy to check that
\begin{equation*}
\begin{split}
&\partial _{x_1 } \tilde {\omega }_1 + \partial _{x_2 } \tilde
{\omega }_2 +\partial _{x_3 } \tilde {\omega }_3 =0\quad \Rightarrow
\quad \partial _{x_1 } \overline {\omega }_1^k +\partial _{x_2 }
\overline
{\omega }_2^k +\partial _{x_3 } \overline {\omega }_3^k =0 \\
&\partial _{x_1 } u_1 +\partial _{x_2 } u_2 +\partial _{x_3 } u_3 =0
\quad\, \Rightarrow \quad \partial _{x_1 } \overline{\overline
u}_1^k +\partial _{x_2 } \overline{\overline u}_2^k +\partial _{x_3 } \overline{\overline u}_3^k =0 \\
\end{split}
\end{equation*}

In the section 3, by means of the Galerkin method and the
compactness imbedding theorem, we can prove the local existence of
the weak solutions of these systems for each $(t_{k-1},\;t_k)$ being
small enough. Below we also interpret $\tilde \omega$ as the
Galerkin approximation of the solution of the problems (5), and
first prove that $\tilde \omega, t \in (0,T)$ belong to
$L^{\infty}(0,T;L^2(\Omega))$. In section 4, an approach of
approximation is used to assert that the solution of (\ref{eq4})
also belongs to $L^{\infty}(0,T;L^2(\Omega))$.
\\

Since
\begin{equation*}
\begin{split}
 &\int_\Omega {(\tilde {\omega }_1 ( \overline{\overline u}_1^k \partial _{x_1 } \overline {\omega
}_1^{k} + \overline{\overline u}_2^k \partial _{x_2 } \overline
\omega _1^{k} + \overline{\overline u}_3^k
\partial _{x_3 } \overline \omega
_1^{k} )} \\
 &\;\,\, +\tilde {\omega }_2 ( \overline{\overline u}_1^k \partial _{x_1 } \overline {\omega }_2^{k} + \overline{\overline u}_2^k
\partial _{x_2 } \overline {\omega }_2^{k} + \overline{\overline u}_3^k \partial _{x_3 } \overline {\omega
}_2^{k} ) \\
 &\;\,\, +\tilde {\omega }_3 ( \overline{\overline u}_1^k \partial _{x_1 } \overline {\omega }_3^{k} + \overline{\overline u}_2^k
\partial _{x_2 } \overline {\omega }_3^{k} + \overline{\overline u}_3^k \partial _{x_3 } \overline {\omega
}_3^{k} )) \\
\end{split}
\end{equation*}
\begin{equation*}
\begin{split}
 &=-\int_\Omega {(\overline {\omega }_1^{k} (\partial _{x_1 } (\tilde {\omega
}_1 \overline{\overline u}_1^k )+\overline {\omega }_1^{k} \partial
_{x_2 } (\tilde {\omega }_1 \overline{\overline u}_2^k
)+\overline {\omega }_1^{k} \partial _{x_3 } (\tilde {\omega }_1 \overline{\overline u}_3^k )} \\
 &\quad \quad \;\;\, +\overline {\omega }_2^{k} \partial _{x_1 } (\tilde {\omega }_2
\overline{\overline u}_1^k )+\overline {\omega }_2^{k} \partial
_{x_2 } (\tilde {\omega }_2 \overline{\overline u}_2^k )+\overline
{\omega }_2^{k} \partial _{x_3 } (\tilde {\omega }_2 \overline{\overline u}_3^k ) \\
 &\quad \quad \;\;\, +\overline {\omega }_3^{k} \partial _{x_1 } (\tilde {\omega }_3
\overline{\overline u}_1^k )+\overline {\omega }_3^{k} \partial
_{x_2 } (\tilde {\omega }_3 \bar u_2^k )+\overline
{\omega }_3^{k} \partial _{x_3 } (\tilde {\omega }_3 \overline{\overline u}_3^k )) \\
 &=-\int_\Omega {(\overline {\omega }_1^{k} \overline{\overline u}_1^k \partial _{x_1 } \tilde {\omega
}_1 +\overline {\omega }_1^{k} \tilde {\omega }_1 \partial _{x_1 }
\overline{\overline u}_1^k +\overline {\omega }_1^{k}
\overline{\overline u}_2^k
\partial _{x_2 } \tilde {\omega }_1 +\overline {\omega }_1^{k}
\tilde {\omega }_1
\partial _{x_2 } \overline{\overline u}_2^k +\overline {\omega }_1^{k} \overline{\overline u}_3^k
\partial _{x_3 } \tilde {\omega }_1 } +\overline {\omega }_1^{k}
\tilde {\omega }_1 \partial _{x_3 } \overline{\overline u}_3^k \\
 &\quad \quad \;\; +\overline {\omega }_2^{k} \overline{\overline u}_1^k \partial _{x_1 } \tilde {\omega }_2
+\overline {\omega }_2^{k} \tilde {\omega }_2 \partial _{x_1 }
\overline{\overline u}_1^k +\overline {\omega }_2^{k}
\overline{\overline u}_2^k
\partial _{x_2 } \tilde {\omega }_2 +\overline {\omega }_2^{k}
\tilde {\omega }_2 \partial _{x_2 } \overline{\overline u}_2^k
+\overline {\omega }_2^k \overline{\overline u}_3^k \partial _{x_3 }
\tilde {\omega }_2 +\overline {\omega }_2^{k} \tilde
{\omega }_2 \partial _{x_3 } \overline{\overline u}_3^k \\
 &\quad \quad \;\; + \overline {\omega }_3^{k} \overline{\overline u}_1^k \partial _{x_1 } \tilde {\omega }_3
+\overline {\omega }_3^{k} \tilde {\omega }_3 \partial _{x_1 }
\overline{\overline u}_1^k +\overline {\omega }_3^{k}
\overline{\overline u}_2^k
\partial _{x_2 } \tilde {\omega }_3 +\overline {\omega }_3^{k}
\tilde {\omega }_3 \partial _{x_2 } \overline{\overline u}_2^k
+\overline {\omega }_3^{k} \overline{\overline u}_3^k \partial _{x_3
} \tilde {\omega }_3 +\overline {\omega }_3^{k} \tilde
{\omega }_3 \partial _{x_3 } \overline{\overline u}_3^k ) \\
 &=-\int_\Omega {(\overline {\omega }_1^{k} \overline{\overline u}_1^k \partial _{x_1 } \tilde {\omega
}_1 +\overline {\omega }_1^{k} \overline{\overline u}_2^k \partial
_{x_2 } \tilde {\omega }_1 +\overline {\omega }_1^{k}
\overline{\overline u}_3^k
\partial _{x_3 } \tilde {\omega }_1 }
\\
 &\quad \quad \;\;\, +\overline {\omega }_2^{k} \overline{\overline u}_1^k \partial _{x_1 } \tilde {\omega }_2
+\overline {\omega }_2^{k} \overline{\overline u}_2^k \partial _{x_2
} \tilde {\omega }_2 +\overline
{\omega }_2^{k} \overline{\overline u}_3^k \partial _{x_3 } \tilde {\omega }_2 \\
 &\quad \quad \;\;\, +\overline {\omega }_3^{k} \overline{\overline u}_1^k \partial _{x_1 } \tilde {\omega }_3
+\overline {\omega }_3^{k} \overline{\overline u}_2^k \partial _{x_2
} \tilde {\omega }_3 +\overline
{\omega }_3^{k} \overline{\overline u}_3^k \partial _{x_3 } \tilde {\omega }_3 ) \\
 \end{split}
\end{equation*}
and
\begin{equation*}
\begin{split}
 &\int_\Omega {(\tilde {\omega }_1 (\overline {\omega }_1^{k} \partial _{x_1 }
\overline{\overline u}_1^k +\overline {\omega }_2^{k} \partial _{x_2
} \overline{\overline u}_1^k +\overline {\omega }_3^{k}
\partial _{x_3 } \overline{\overline u}_1^k } ) \\
 &\;\;+\tilde {\omega }_2 (\overline {\omega }_1^{k} \partial _{x_1 } \overline{\overline u}_2^k +\overline
{\omega }_2^{k} \partial _{x_2 } \overline{\overline u}_2^k
+\overline {\omega }_3^{k}
\partial
_{x_3 } \overline{\overline u}_2^k ) \\
 &\;\;+\tilde {\omega }_3 (\overline {\omega }_1^{k} \partial _{x_1 } \overline{\overline u}_3^k +\overline
{\omega }_2^{k} \partial _{x_2 } \overline{\overline u}_3^k
+\overline {\omega }_3^{k}
\partial
_{x_3 } \overline{\overline u}_3^k )) \\
 &=-\int_\Omega {\,(\overline{\overline u}_1^k \partial _{x_1 } (\tilde {\omega }_1 \overline {\omega
}_1^{k} )+ \overline{\overline u}_1^k \partial _{x_2 } (\tilde
{\omega }_1 \overline {\omega }_2^{k}
)+ \overline{\overline u}_1^k \partial _{x_3 } (\tilde {\omega }_1 \overline {\omega }_3^{k} )} \\
 &\quad \quad \;\;\; + \overline{\overline u}_2^k \partial _{x_1 } (\tilde {\omega }_2 \overline {\omega
}_1^{k} )+ \overline{\overline u}_2^k \partial _{x_2 } (\tilde
{\omega }_2 \overline {\omega }_2^{k}
)+ \overline{\overline u}_2^k \partial _{x_3 } (\tilde {\omega }_2 \overline {\omega }_3^{k} ) \\
 &\quad \quad \;\;\; + \overline{\overline u}_3^k \partial _{x_1 } (\tilde {\omega }_3 \overline {\omega
}_1^{k} )+ \overline{\overline u}_3^k \partial _{x_2 } (\tilde
{\omega }_3 \overline {\omega }_2^{k}
)+ \overline{\overline u}_3^k \partial _{x_3 } (\tilde {\omega }_3 \overline {\omega }_3^{k} )) \\
 &=-\int_\Omega {(\overline {\omega }_1^{k} \overline{\overline u}_1^k \partial _{x_1 } \tilde {\omega
}_1 +\tilde {\omega }_1 \overline{\overline u}_1^k \partial _{x_1 }
\overline {\omega }_1^{k} +\overline {\omega }_2^{k}
\overline{\overline u}_1^k
\partial _{x_2 } \tilde {\omega }_1 +\tilde {\omega }_1 \overline{\overline u}_1^k
\partial _{x_2 } \overline {\omega }_2^{k} +\overline {\omega
}_3^{k} \overline{\overline u}_1^k
\partial _{x_3 } \tilde {\omega }_1 }
+\tilde {\omega }_1 \overline{\overline u}_1^k \partial _{x_3 } \overline {\omega }_3^{k} \\
 &\quad \quad \;\;\, +\overline {\omega }_1^{k} \overline{\overline u}_2^k \partial _{x_1 } \tilde {\omega }_2
+\tilde {\omega }_2 \overline{\overline u}_2^k \partial _{x_1 }
\overline {\omega }_1^{k} +\overline {\omega }_2^{k}
\overline{\overline u}_2^k \partial _{x_2 } \tilde {\omega }_2
+\tilde {\omega }_2 \overline{\overline u}_2^k \partial _{x_2 }
\overline {\omega }_2^{k} +\overline {\omega }_3^{k}
\overline{\overline u}_2^k
\partial _{x_3 } \tilde {\omega }_2 +\tilde {\omega }_2 \overline{\overline u}_2^k \partial _{x_3 }
\overline {\omega }_3^{k} \\
 &\quad \quad \;\;\, +\overline {\omega }_1^{k} \overline{\overline u}_3^k \partial _{x_1 } \tilde {\omega }_3
+\tilde {\omega }_3 \overline{\overline u}_3^k \partial _{x_1 }
\overline {\omega }_1^{k} +\overline {\omega }_2^{k}
\overline{\overline u}_3^k
\partial _{x_2 } \tilde {\omega }_3 +\tilde {\omega }_3 \overline{\overline u}_3^k
\partial _{x_2 } \overline {\omega }_2^{k} +\overline {\omega
}_3^{k} \overline{\overline u}_3^k
\partial _{x_3 } \tilde {\omega }_3 +\tilde {\omega }_3 \overline{\overline u}_3^k
\partial _{x_3 } \overline {\omega }_3^{k} ) \\
 &=-\int_\Omega {(\overline {\omega }_1^{k} \overline{\overline u}_1^k \partial _{x_1 } \tilde {\omega
}_1 +\overline {\omega }_2^{k} \overline{\overline u}_1^k \partial
_{x_2 } \tilde {\omega }_1
+\overline {\omega }_3^{k} \overline{\overline u}_1^k \partial _{x_3 } \tilde {\omega }_1 } \\
 &\quad \quad \;\;\, +\overline {\omega }_1^{k} \overline{\overline u}_2^k \partial _{x_1 } \tilde {\omega }_2
+\overline {\omega }_2^{k} \overline{\overline u}_2^k \partial _{x_2
} \tilde {\omega }_2 +\overline
{\omega }_3^{k} \overline{\overline u}_2^k \partial _{x_3 } \tilde {\omega }_2 \\
 &\quad \quad \;\;\, + \overline {\omega }_1^{k} \overline{\overline u}_3^k \partial _{x_1 } \tilde {\omega }_3
+\overline {\omega }_2^{k} \overline{\overline u}_3^k \partial _{x_2
} \tilde {\omega }_3 +\overline
{\omega }_3^{k} \overline{\overline u}_3^k \partial _{x_3 } \tilde {\omega }_3 ) \\
 \end{split}
\end{equation*}
\\
\[
\int_\Omega {(\tilde {\omega }_1 \partial _{x_1 } q+\tilde {\omega }_2
\partial _{x_2 } q+\tilde {\omega }_3 \partial _{x_3 } q)} =-\int_\Omega
{q\,(\partial _{x_1 } \tilde {\omega }_1 +\partial _{x_2 } \tilde {\omega
}_2 +\partial _{x_3 } \tilde {\omega }_3 )} =0
\]
\\
furthermore
\[
 \int_\Omega {\tilde {\omega }_i \Delta \tilde {\omega }_i } =\int_\Omega
{\tilde {\omega }_i (\partial _{x_1 }^2 \tilde {\omega }_i +\partial
_{x_2 }^2 \tilde {\omega }_i +\partial _{x_3 }^2 \tilde {\omega }_i
)} =-\int_\Omega {((\partial _{x_1 } \tilde {\omega }_i
)^2+(\partial _{x_2 } \tilde {\omega }_i )^2+(\partial _{x_3 }
\tilde {\omega }_i )^2)}
\]
\\
Then from (\ref{eq5}) we have
\begin{equation*}
\begin{split}
 &\int_\Omega {\tilde {\omega }_1 \partial _t \tilde {\omega }_1 }
\;\,+\int_\Omega {\tilde {\omega }_1 (\overline{\overline u}_1^k
\partial _{x_1 } \overline {\omega }_1^{k} + \overline{\overline u}_2^k \partial
_{x_2 } \overline {\omega }_1^{k} + \overline{\overline u}_3^k
\partial _{x_3 }
\overline {\omega }_1^{k} )} \\
 &\quad \quad \quad \quad \quad \quad \quad \quad \quad -\int_\Omega {\tilde
{\omega }_1 (\overline {\omega }_1^{k} \partial _{x_1 }
\overline{\overline u}_1^k +\overline {\omega }_2^{k} \partial _{x_2
} \overline{\overline u}_1^k +\overline {\omega }_3^{k}
\partial _{x_3 } \overline{\overline u}_1^k )} \,+ \int_\Omega {\tilde {\omega }_1 \partial _{x_1 }
q} =\int_\Omega {\tilde {\omega }_1
\Delta \tilde {\omega }_1 } \\
 &\int_\Omega {\tilde {\omega }_2 \partial _t \tilde {\omega }_2 }
+\int_\Omega {\tilde {\omega }_2 (\overline{\overline u}_1^k
\partial _{x_1 } \overline {\omega }_2^{k} + \overline{\overline u}_2^k \partial
_{x_2 } \overline {\omega }_2^{k} + \overline{\overline u}_3^k
\partial _{x_3 }
\overline {\omega }_2^{k} )} \\
 &\quad \quad \quad \quad \quad \quad \quad \quad \quad -\int_\Omega {\tilde
{\omega }_2 (\overline {\omega }_1^{k} \partial _{x_1 }
\overline{\overline u}_2^k +\overline {\omega }_2^{k} \partial _{x_2
} \overline{\overline u}_2^k +\overline {\omega }_3^{k}
\partial _{x_3 } \overline{\overline u}_2^k )} + \int_\Omega {\tilde {\omega }_2 \partial _{x_2 }
q} =\int_\Omega {\tilde {\omega }_2
\Delta \tilde {\omega }_2 } \\
 &\int_\Omega {\tilde {\omega }_3 \partial _t \tilde {\omega }_3 }
+\int_\Omega {\tilde {\omega }_3 (\overline{\overline u}_1^k
\partial _{x_1 } \overline {\omega }_3^{k} + \overline{\overline u}_2^k \partial
_{x_2 } \overline {\omega }_3^{k} + \overline{\overline u}_3^k
\partial _{x_3 }
\overline {\omega }_3^{k} )} \\
 &\quad \quad \quad \quad \quad \quad \quad \quad \quad -\int_\Omega {\tilde
{\omega }_3 (\overline {\omega }_1^{k} \partial _{x_1 }
\overline{\overline u}_3^k +\overline {\omega }_2^{k} \partial _{x_2
} \overline{\overline u}_3^k + \overline {\omega }_3^{k}
\partial _{x_3 } \overline{\overline u}_3^k )} + \int_\Omega { \tilde {\omega }_3 \partial _{x_3 }
q} =\int_\Omega {\tilde {\omega }_3
\Delta \tilde {\omega }_3 } \\
 \end{split}
\end{equation*}
so that
\begin{equation*}
\begin{split}
 &\frac{1}{2}\partial _t \int_\Omega {(\tilde {\omega }_1^2 +\tilde {\omega
}_2^2 +\tilde {\omega }_3^2 )\;} +\,\,\,\int_\Omega {((\partial _{x_1 }
\tilde {\omega }_1 )^2+(\partial _{x_2 } \tilde {\omega }_1 )^2+(\partial
_{x_3 } \tilde {\omega }_1 )^2)} \\
 &\quad \quad \quad \quad \quad \quad \quad \quad \quad \quad \qquad +(\partial
_{x_1 } \tilde {\omega }_2 )^2+(\partial _{x_2 } \tilde {\omega }_2
)^2+(\partial _{x_3 } \tilde {\omega }_2 )^2 \\
 &\quad \quad \quad \quad \quad \quad \quad \quad \quad \quad \qquad +(\partial
_{x_1 } \tilde {\omega }_3 )^2+(\partial _{x_2 } \tilde {\omega }_3
)^2+(\partial _{x_3 } \tilde {\omega }_3 )^2) \\
 &\quad \quad -\int_\Omega {(\overline {\omega }_1^{k} \overline{\overline u}_1^k \partial _{x_1 }
\tilde {\omega }_1 +\overline {\omega }_1^{k} \overline{\overline
u}_2^k \partial _{x_2 } \tilde {\omega }_1 +\overline {\omega
}_1^{k} \overline{\overline u}_3^k
\partial _{x_3 } \tilde {\omega
}_1 } \\
 &\quad \quad \quad \;\; +\overline {\omega }_2^{k} \overline{\overline u}_1^k \partial _{x_1 } \tilde
{\omega }_2 +\overline {\omega }_2^{k} \overline{\overline u}_2^k
\partial _{x_2 } \tilde {\omega }_2
+\overline {\omega }_2^{k} \overline{\overline u}_3^k \partial _{x_3 } \tilde {\omega }_2 \\
 &\quad \quad \quad \;\; +\overline {\omega }_3^{k} \overline{\overline u}_1^k \partial _{x_1 } \tilde
{\omega }_3 +\overline {\omega }_3^{k} \overline{\overline u}_2^k
\partial _{x_2 } \tilde {\omega
}_3 \,\,+\overline {\omega }_3^{k} \overline{\overline u}_3^k \partial _{x_3 } \tilde {\omega }_3 ) \\
 &\quad \quad +\int_\Omega {(\overline {\omega }_1^{k} \overline{\overline u}_1^k \partial _{x_1 }
\tilde {\omega }_1 +\overline {\omega }_2^{k} \overline{\overline
u}_1^k \partial _{x_2 } \tilde {\omega }_1 +\overline {\omega
}_3^{k} \overline{\overline u}_1^k
\partial _{x_3 } \tilde {\omega
}_1 } \\
 &\quad \quad \quad \;\;\, +\overline {\omega }_1^{k} \overline{\overline u}_2^k \partial _{x_1 } \tilde
{\omega }_2 +\overline {\omega }_2^{k} \overline{\overline u}_2^k
\partial _{x_2 } \tilde {\omega }_2
+\overline {\omega }_3^{k} \overline{\overline u}_2^k \partial _{x_3 } \tilde {\omega }_2 \\
 &\quad \quad \quad \;\;\, +\overline {\omega }_1^{k} \overline{\overline u}_3^k \partial _{x_1 } \tilde
{\omega }_3 +\overline {\omega }_2^{k} \overline{\overline u}_3^k
\partial _{x_2 } \tilde {\omega }_3 +\overline {\omega }_3^{k} \overline{\overline u}_3^k \partial _{x_3 } \tilde {\omega }_3 )=0
\\
 \end{split}
\end{equation*}
it follows that
\begin{equation*}
\begin{split}
 &\partial _t \int_\Omega {(\tilde {\omega }_1^2 +\tilde {\omega }_2^2 +\tilde {\omega
}_3^2 )} \;\,+\,\;2\;\; {\int_\Omega {\;((\partial _{x_1 } \tilde
{\omega }_1 )^2 +(\partial _{x_2 } \tilde {\omega }_1
)^2+(\partial _{x_3 } \tilde {\omega }_1 )^2} } \\
 &\qquad \qquad \qquad \quad \quad \quad \quad \quad \quad \quad
\;\; +(\partial _{x_1 } \tilde {\omega }_2 )^2+(\partial _{x_2 }
\tilde
{\omega }_2 )^2+(\partial _{x_3 } \tilde {\omega }_2 )^2 \\
 &\qquad \qquad \qquad \quad \quad \quad \quad \quad \quad \quad
\;\; +(\partial _{x_1 } \tilde {\omega }_3 )^2+(\partial _{x_2 }
\tilde
{\omega }_3 )^2+(\partial _{x_3 } \tilde {\omega }_3 )^2) \\
 &\qquad \qquad \qquad \qquad \le \; 2 \; {\int_\Omega
{(\overline {\omega }_1^{k^2} \overline{\overline u}_1^{k^2}
+\overline {\omega }_1^{k^2} \overline{\overline u}_2^{k^2}
+\overline {\omega
}_1^{k^2} \overline{\overline u}_3^{k^2} } } \\
 &\qquad \qquad \quad \quad \,\quad \;\; \qquad
\quad \quad + \overline {\omega }_2^{k^2} \overline{\overline
u}_1^{k^2} +\overline {\omega }_2^{k^2} \overline{\overline
u}_2^{k^2}
+\overline {\omega }_2^{k^2} \overline{\overline u}_3^{k^2} \\
 &\qquad \qquad \qquad \qquad \;\; \quad
\quad \quad +\overline {\omega }_3^{k^2} \overline{\overline
u}_1^{k^2} \,+\overline {\omega }_3^{k^2}
\overline{\overline u}_2^{k^2} \,+ \overline {\omega }_3^{k^2} \overline{\overline u}_3^{k^2} ) \\
\end{split}
\end{equation*}
\begin{equation*}
\begin{split}
 &\quad \qquad \qquad \qquad
\;\;\;\,+\;\,2 \; {\int_\Omega {(\overline {\omega }_1^{k^2}
\overline{\overline u}_1^{k^2} +\overline {\omega }_2^{k^2} \overline{\overline u}_1^{k^2} + \overline {\omega }_3^{k^2} \overline{\overline u}_1^{k^2} } } \\
 &\qquad \qquad \qquad \qquad \quad \;\;
\quad \quad +\overline {\omega }_1^{k^2} \overline{\overline
u}_2^{k^2} + \overline {\omega }_2^{k^2} \overline{\overline
u}_2^{k^2}
+\overline {\omega }_3^{k^2} \overline{\overline u}_2^{k^2} \\
 &\qquad \qquad \qquad \qquad \quad \;\;
\quad \quad + \overline {\omega }_1^{k^2} \overline{\overline
u}_3^{k^2} \,+ \overline {\omega }_2^{k^2}
\overline{\overline u}_3^{k^2} \,+ \overline {\omega }_3^{k^2} \overline{\overline u}_3^{k^2} ) \\
 &\quad \quad \quad \quad \quad \quad \qquad  + \; {\int_\Omega
{((\partial _{x_1 } \tilde {\omega }_1 )^2+(\partial _{x_2 } \tilde {\omega
}_1 )^2+(\partial _{x_3 } \tilde {\omega }_1 )^2} } \\
\end{split}
\end{equation*}
\begin{equation*}
\begin{split}
 &\qquad \quad \quad \quad \quad \quad \quad \quad \quad \;\;\; +(\partial _{x_1
} \tilde {\omega }_2 )^2+(\partial _{x_2 } \tilde {\omega }_2 )^2+(\partial
_{x_3 } \tilde {\omega }_2 )^2 \\
 &\qquad \quad \quad \quad \quad \quad \quad \quad \quad \;\;\; +(\partial _{x_1
} \tilde {\omega }_3 )^2+(\partial _{x_2 } \tilde {\omega }_3 )^2+(\partial
_{x_3 } \tilde {\omega }_3 )^2) \\
 \end{split}
\end{equation*}
by using Young inequality: $uv\le \frac{1}{4}u^2+v^2$.
\\

Thus,
\begin{equation}
\label{eq6}
\begin{split}
 &\int_\Omega {(\tilde {\omega }_1^2 +\tilde {\omega }_2^2 +\tilde {\omega
}_3^2 )} +\int_{t_{k-1} }^t {\int_\Omega {((\partial _{x_1 } \tilde {\omega
}_1 )^2+(\partial _{x_2 } \tilde {\omega }_1 )^2+(\partial _{x_3 } \tilde
{\omega }_1 )^2} } \\
 &\quad \qquad \qquad \qquad \qquad \qquad \quad \;\;\, +(\partial
_{x_1 } \tilde {\omega }_2 )^2+(\partial _{x_2 } \tilde {\omega }_2
)^2+(\partial _{x_3 } \tilde {\omega }_2 )^2 \\
 &\quad \qquad \qquad \qquad \qquad \qquad \quad \;\;\, +(\partial _{x_1
} \tilde {\omega }_3 )^2+(\partial _{x_2 } \tilde {\omega }_3 )^2+(\partial
_{x_3 } \tilde {\omega }_3 )^2) \\
 &\qquad \le \int_\Omega {(\tilde {\omega }_1^{k-1^2} +\tilde {\omega }_2^{k-1^2} +\tilde
{\omega }_3^{k-1^2} )} \;\;\; + \\
 &\qquad \qquad + 4\, \int_{t_{k-1} }^t  (\,\left\| {\overline{\overline u}_1^k } \right\|_{L^{\infty}(\Omega )}^2 +\left\| {\overline{\overline u}_2^k }
\right\|_{L^{\infty}(\Omega )}^2 +\left\| {\overline{\overline
u}_3^k } \right\|_{L^{\infty}(\Omega )}^2 ) \;\;  \int_\Omega {(\,
\overline {\omega }_1^{k^2} +\overline {\omega }_2^{k^2} +\overline
{\omega }_3^{k^2} )} \\
 &\qquad \le \int_\Omega {(\tilde {\omega }_1^{k-1^2} +\tilde {\omega }_2^{k-1^2} +\tilde
{\omega }_3^{k-1^2} )} \;\;\; + \\
 &\qquad \qquad + 4\; \Delta t_k \; (\,\left\| {\overline{\overline u}_1^k } \right\|_{L^{\infty}(\Omega )}^2 +\left\| {\overline{\overline u}_2^k }
\right\|_{L^{\infty}(\Omega )}^2 +\left\| {\overline{\overline
u}_3^k } \right\|_{L^{\infty}(\Omega )}^2 ) \; (\,\left\| {\overline
{\omega }_1^{k} } \right\|_{L^2 (\Omega )}^2 +\left\| {\overline
{\omega }_2^{k} } \right\|_{L^2 (\Omega )}^2 +\left\| {\overline
{\omega
}_3^{k} } \right\|_{L^2 (\Omega )}^2 ) \\
 \end{split}
\end{equation}

Note that
\begin{equation*}
\begin{split}
 &\left\| {\overline {\omega }_i^{k} } \right\|_{L^2(\Omega )}^2 =\int_\Omega
{\left( {\frac{1}{\Delta t_{k} }\int_{t_{k-1} }^{t_{k} } {\tilde
{\omega }_i (x,t)dt} } \right)} ^2\le \frac{1}{\Delta t_{k}^2
}\int_\Omega {\Delta
t_{k} \int_{t_{k-1} }^{t_{k} } {\tilde {\omega }_i^2 (x,t)dt} } \\
 &\quad \quad \quad \quad \;\; =\frac{1}{\Delta t_{k} }\int_{t_{k-1} }^{t_{k}
} {\left\| {\tilde {\omega }_i } \right\|_{L^2(\Omega )}^2 } \\
 \end{split}
\end{equation*}
and similarly
\[
 \left\| {\overline u_i^{k} } \right\|_{L^2(\Omega )}^2
\le \frac{1}{\Delta t_{k} }\int_{t_{k-1} }^{t_{k}
} {\left\| u_i \right\|_{L^2(\Omega )}^2 }, \qquad i=1,2,3 \\
\]
In addition, a convolution inequality in [1] is applied to get
\begin{equation*}
\begin{split}
 &\left\| {\overline{\overline u} _i^k } \right\|_{L^\infty (\Omega )}^2 =
\left\| {J_\varepsilon * \overline u _i^k (x)} \right\|_{L^\infty
(\Omega
)}^2 \\
 &\qquad \qquad \le \left\| {J_\varepsilon } \right\|_{L^2(B_\varepsilon )}^2
\;\left\| {\overline u _i^k } \right\|_{L^2(\Omega )}^2 \;\; \le \;\frac{1}{\mu_\varepsilon}
\;\mathop {\sup }\limits_{(t_{k - 1} ,t_k )} \left\| {u_i }
\right\|_{L^2(\Omega )}^2
\\
\end{split}
\end{equation*}
where $B_\varepsilon = \{x:\left| x \right| < \varepsilon \}$ and
$\left\| {\overline{\overline u} _i^k } \right\|_{L^\infty (\Omega
)} = \left\| {\overline{\overline u} _i^k } \right\|_{L^\infty
(\mathbb R^3)} $, the quantity $\mu_\varepsilon \to 0$ as $\varepsilon \to 0$. We still need further assuming that $\varepsilon \to 0$
 and $\frac{\Delta t_k}{\mu_\varepsilon} \to 0$ as $k \to \infty$ or $\Delta t_k \to 0$.
\\

From (\ref{eq6}) we have
\begin{equation*}
\begin{split}
 &\int_\Omega {(\tilde {\omega }_1^2 + \tilde {\omega }_2^2 +\tilde {\omega
}_3^2 )} \;+\; \int_{t_{k-1} }^t {(\,\left\| {\nabla \tilde {\omega
}_1 } \right\|_{L^2(\Omega )}^2 +\left\| {\nabla \tilde {\omega }_2
} \right\|_{L^2(\Omega )}^2 +\left\| {\nabla \tilde {\omega }_3 }
\right\|_{L^2(\Omega )}^2 )} \\
 &\qquad  \le \int_\Omega {(\tilde {\omega }_1^{k-1^2} +\tilde {\omega }_2^{k-1^2} +\tilde
{\omega }_3^{k-1^2} )} \;\; + \\
 &\qquad \qquad  + \;\frac{4\Delta t_{k} }{\mu_\varepsilon} \; \mathop {\sup }\limits_{(t_{k-1} ,t_k )} \left\{
{\left\| {u_1 } \right\|_{L^2(\Omega )}^2 +\left\| {u_2 }
\right\|_{L^2(\Omega )}^2 +\left\| {u_3 } \right\|_{L^2(\Omega )}^2
} \right\}  \cdot \mathop {\sup }\limits_{(t_{k-1} ,t )}\; \int_\Omega {(\tilde {\omega }_1^2 + \tilde {\omega }_2^2 +\tilde {\omega
}_3^2 )}  \\
 \end{split}
\end{equation*}

By (\ref{eq2}) we have
\begin{equation*}
\begin{split}
  &\mathop {\sup }\limits_{(t_{k-1} ,t_k )} \left\{
{\left\| {u_1 } \right\|_{L^2(\Omega )}^2 +\left\| {u_2 }
\right\|_{L^2(\Omega )}^2 +\left\| {u_3 } \right\|_{L^2(\Omega )}^2
} \right\}  \\
  &\qquad  \le K_0 = \mathop {\sup }\limits_{t\in (0,T)} \; \int_\Omega {(u_1^2
+u_2^2 +u_3^2 )} + \int_0
^T {\left( {\left\| {\nabla u_1 }
\right\|_{L^2(\Omega )}^2 +\left\| {\nabla u_2 }
\right\|_{L^2(\Omega )}^2 +\left\| {\nabla u_3 }
\right\|_{L^2(\Omega )}^2 } \right)}\; < +\infty   \\
 \end{split}
\end{equation*}
Thus,
\begin{equation*}
\begin{split}
 & {\left( 1-4K_0 \frac{\Delta t_{k} }{\mu_\varepsilon} \right)}\; \mathop {\sup }\limits_{t\in (t_{k-1} ,t_k )} \int_\Omega {(\tilde {\omega
}_1^2 \;+ \tilde {\omega }_2^2 +\tilde {\omega }_3^2 )} \;\; + \\
 &\quad + \int_{t_{k-1}
}^{t_k } {\left( {\left\| {\nabla \tilde {\omega }_1 }
\right\|_{L^2(\Omega )}^2 +\left\| {\nabla \tilde {\omega }_2 }
\right\|_{L^2(\Omega )}^2 +\left\| {\nabla \tilde {\omega }_3 }
\right\|_{L^2(\Omega )}^2 } \right)} \le \int_\Omega {(\tilde {\omega }_1^{k-1^2} +\tilde {\omega }_2^{k-1^2}
+\tilde {\omega }_3^{k-1^2} )}  \\
 \end{split}
\end{equation*}

Now we set
\begin{equation*}
\begin{split}
 & M_0 =\int_\Omega {(\omega _{10}^2 +\omega _{20}^2 +\omega _{30}^2 )}  \\
 & M_k = \mathop {\sup }\limits_{t\in (t_{k-1} ,t_k ) } \;\int_\Omega {(\tilde
{\omega }_1^2 +\tilde {\omega }_2^2 +\tilde {\omega }_3^2 )} \\
 & \delta_k = \int_{t_{k-1}
}^{t_k } {\left( {\left\| {\nabla \tilde {\omega }_1 }
\right\|_{L^2(\Omega )}^2 +\left\| {\nabla \tilde {\omega }_2 }
\right\|_{L^2(\Omega )}^2 +\left\| {\nabla \tilde {\omega }_3 }
\right\|_{L^2(\Omega )}^2 } \right)}  \\
 &\qquad \qquad \qquad \qquad \qquad \qquad \qquad   k=1,\cdots ,N \\
\end{split}
\end{equation*}
then we have
\[
{\left( 1-4K_0 \frac{\Delta t_{k} }{\mu_\varepsilon} \right)}\;M_k\; + \;\delta_k \le\; M_{k-1}
\]
\\

The partition is assumed to be fine enough. Because of the local existence of Galerkin solution in section 3 and the absolute continuity of integration with respect to $t$, it is valid that $\delta_k \to 0$
as $\Delta t_k \to 0$.
\\

We may first consider the case that
\[
{M_{k-1}} \frac{\Delta t_k }{\delta_k} \to 0, \quad \mbox{as} \;\;\Delta t_k \to 0
\]
which may be a subsequence $k'$, still denoted $k$. At this time, we can choose $\varepsilon_k$ on each $(t_{k-1},\;t_k)$ such that
\[
\mu_{\varepsilon_k} = 4K_0\,M_{k-1}\;\frac{\Delta t_k }{\delta_k} \quad \mbox{and} \quad 1-4K_0\,\frac{\Delta t_k }{\mu_{\varepsilon_k}} \ge \frac{1}{2}
\]
$\varepsilon = \mathop{\max }\limits_k \{\varepsilon_k \}$.
\\

Then we obtain
\[
{\left( 1-4K_0 \frac{\Delta t_k }{\mu_{\varepsilon_k}} \right)} \;M_k\; + \;\delta_k \;= {\left( 1- \frac{\delta_k }{M_{k-1}} \right)} \;M_k\; + \;\delta_k\; \le\; M_{k-1}
\]
it follows that $M_k \le M_{k-1}$.
\\
\\

Otherwise, $\delta_k \le O(\Delta t_k) {M_{k-1}}$. This leads to that
\begin{equation*}
\begin{split}
 &\{\;\left\| {\nabla {\overline {\omega }}_1^k }
\right\|_{L^2(\Omega )}^2 + \left\| {\nabla {\overline
{\omega }}_2^k } \right\|_{L^2(\Omega )}^2 + \left\| {\nabla {\overline {\omega }}_3^k }
\right\|_{L^2(\Omega )}^2 \}  \\
 &\qquad  \le \frac{1}{\Delta t_k} \int_{t_{k-1}
}^{t_k } {\left( {\left\| {\nabla \tilde {\omega }_1 }
\right\|_{L^2(\Omega )}^2 +\left\| {\nabla \tilde {\omega }_2 }
\right\|_{L^2(\Omega )}^2 +\left\| {\nabla \tilde {\omega }_3 }
\right\|_{L^2(\Omega )}^2 } \right)} \; \le O(1) {M_{k-1}} \\
\end{split}
\end{equation*}
Since these $(t_{k-1},\;t_k)$ are of finite length, the number of them is finite. According to Cauchy-Schwartz inequality, similar to (6), we have
\\
\begin{equation*}
\begin{split}
 &\partial _t \int_\Omega {(\tilde {\omega }_1^2 + \tilde {\omega }_2^2 +
\tilde {\omega }_3^2 )} \;\; + \int_\Omega {[\,(\partial _{x_1 }
\tilde {\omega }_1 )^2 + (\partial _{x_2 } \tilde {\omega }_1 )^2 +
(\partial _{x_3
} \tilde {\omega }_1 )^2} \\
 &\quad \quad \quad \quad \quad \quad \quad \qquad \qquad  + (\partial
_{x_1 } \tilde {\omega }_2 )^2 + (\partial _{x_2 } \tilde {\omega
}_2 )^2 +
(\partial _{x_3 } \tilde {\omega }_2 )^2 \\
 &\quad \quad \quad \quad \quad \quad \quad \qquad \qquad  + (\partial
_{x_1 } \tilde {\omega }_3 )^2 + (\partial _{x_2 } \tilde {\omega
}_3 )^2 +
(\partial _{x_3 } \tilde {\omega }_3 )^2\,] \\
 &\le 4\,\{\,(\int_\Omega {\overline {\overline {u}}_1^{k^4} } )^{\frac{1}{2}}(\int_\Omega {
{\overline {\omega }}_1^{k^4} } )^{\frac{1}{2}} + (\int_\Omega
{\overline {\overline {u}}_1^{k^4} } )^{\frac{1}{2}}(\int_\Omega {
{\overline {\omega }}_2^{k^4} } )^{\frac{1}{2}} + (\int_\Omega
{\overline {\overline {u}}_1^{k^4} } )^{\frac{1}{2}}(\int_\Omega {
{\overline {\omega }}_3^{k^4} } )^{\frac{1}{2}}
\\
 &\qquad \;\;\; (\int_\Omega {\overline {\overline {u}}_2^{k^4} } )^{\frac{1}{2}}(\int_\Omega
{ {\overline {\omega }}_1^{k^4} } )^{\frac{1}{2}} +
(\int_\Omega {\overline {\overline {u}}_2^{k^4} } )^{\frac{1}{2}}(\int_\Omega
{ {\overline {\omega }}_2^{k^4} } )^{\frac{1}{2}} +
(\int_\Omega {\overline {\overline {u}}_2^{k^4} } )^{\frac{1}{2}}(\int_\Omega
{ {\overline {\omega }}_3^{k^4} } )^{\frac{1}{2}}
\\
 &\qquad \;\;\; (\int_\Omega {\overline {\overline {u}}_3^{k^4} } )^{\frac{1}{2}}(\int_\Omega
{ {\overline {\omega }}_1^{k^4} } )^{\frac{1}{2}} +
(\int_\Omega {\overline {\overline {u}}_3^{k^4} } )^{\frac{1}{2}}(\int_\Omega
{ {\overline {\omega }}_2^{k^4} } )^{\frac{1}{2}} +
(\int_\Omega {\overline {\overline {u}}_3^{k^4} } )^{\frac{1}{2}}(\int_\Omega
{ {\overline {\omega }}_3^{k^4} }
)^{\frac{1}{2}}\,\} \\
 &= 4\,\{\,\,\left\| {\overline {\overline {u}}_1^k } \right\|_{L^4(\Omega )}^2 (\,\left\|
{ {\overline {\omega }}_1^k } \right\|_{L^4(\Omega )}^2 +
\left\| { {\overline {\omega }}_2^k } \right\|_{L^4(\Omega
)}^2 + \left\| { {\overline {\omega
}}_3^k } \right\|_{L^4(\Omega )}^2 ) \\
 &\quad \;\; + \left\| {\overline {\overline {u}}_2^k } \right\|_{L^4(\Omega )}^2 (\,\left\|
{ {\overline {\omega }}_1^k } \right\|_{L^4(\Omega )}^2 +
\left\| { {\overline {\omega }}_2^k } \right\|_{L^4(\Omega
)}^2 + \left\| { {\overline {\omega
}}_3^k } \right\|_{L^4(\Omega )}^2 ) \\
 \end{split}
\end{equation*}
\begin{equation*}
\begin{split}
 &\quad \;\; + \left\| {\overline {\overline {u}
 }_3^k } \right\|_{L^4(\Omega )}^2 (\,\left\|
{ {\overline {\omega }}_1^k } \right\|_{L^4(\Omega )}^2 +
\left\| { {\overline {\omega }}_2^k } \right\|_{L^4(\Omega
)}^2 + \left\| { {\overline {\omega
}}_3^k } \right\|_{L^4(\Omega )}^2 )\,\} \\
 &= 4\,(\;\left\| {\overline {\overline {u}}_1^k } \right\|_{L^4(\Omega )}^2 + \left\| {\overline {\overline
{u}}_2^k } \right\|_{L^4(\Omega )}^2 + \left\| {\overline {\overline {u}}_3^k }
\right\|_{L^4(\Omega )}^2 )\; (\;\left\| { {\overline
{\omega }}_1^k } \right\|_{L^4(\Omega )}^2 + \left\| {
{\overline {\omega }}_2^k } \right\|_{L^4(\Omega )}^2 + \left\|
{ {\overline {\omega }}_3^k }
\right\|_{L^4(\Omega )}^2 ) \\
 \end{split}
\end{equation*}
\\

From Sobolev imbedding theorem in [1], there exists a constant $C_1
> 0$ independent of $\omega $ such that
\[
\left\| { {\overline {\omega }}_i^k } \right\|_{L^4(\Omega
)}^2 \le C_1 \{\;\left\| { {\overline {\omega }}_i^k }
\right\|_{L^2(\Omega )}^2 + \left\| {\nabla  {\overline
{\omega }}_i^k } \right\|_{L^2(\Omega )}^2 \},\quad \quad i = 1,2,3
\]
Therefore,
\begin{equation*}
\begin{split}
 &\int_\Omega {(\tilde {\omega }_1^2 + \tilde {\omega }_2^2 + \tilde {\omega
}_3^2 )} \;\; + \int_{t_{k - 1} }^t {(\;\left\| {\nabla \tilde
{\omega }_1 } \right\|_{L^2(\Omega )}^2 + \left\| {\nabla \tilde
{\omega }_2 } \right\|_{L^2(\Omega )}^2 + \left\| {\nabla \tilde
{\omega }_3 }
\right\|_{L^2(\Omega )}^2 )} \\
 &\le \int_\Omega {(\tilde {\omega }_1^{k - 1^2} + \tilde {\omega }_2^{k -
1^2} + \tilde {\omega }_3^{k - 1^2} )} \;\; + \\
 &+ \,C_2\,(\;\left\| {\overline{\overline {u}}_1^k } \right\|_{L^2(\Omega )}^2 + \left\|
{\overline{\overline {u}}_2^k } \right\|_{L^2(\Omega )}^2 + \left\| {\overline{\overline
{u}}_3^k } \right\|_{L^2(\Omega )}^2 + \left\| {\nabla \overline{\overline
{u}}_1^k } \right\|_{L^2(\Omega )}^2 + \left\| {\nabla \overline{\overline
{u}}_2^k } \right\|_{L^2(\Omega )}^2 + \left\| {\nabla \overline{\overline
{u}}_3^k }
\right\|_{L^2(\Omega )}^2 ) \\
 &\quad \quad \times \int_{t_{k - 1} }^t {\int_\Omega {( {\overline {\omega
}}_1^{k^2} +  {\overline {\omega }}_2^{k^2} +
{\overline {\omega }}_3^{k^2}
)} } \;\; + \\
 &+ \,C_2\,\int_{t_{k - 1} }^t {(\;\left\| {\overline{\overline {u}}_1^k }
\right\|_{L^2(\Omega )}^2 + \left\| {\overline{\overline {u}}_2^k }
\right\|_{L^2(\Omega )}^2 + } \left\| {\overline{\overline {u}}_3^k }
\right\|_{L^2(\Omega )}^2 + \left\| {\nabla \overline{\overline {u}}_1^k }
\right\|_{L^2(\Omega )}^2 + \left\| {\nabla \overline{\overline {u}}_2^k }
\right\|_{L^2(\Omega )}^2 + \left\| {\nabla \overline{\overline {u}}_3^k }
\right\|_{L^2(\Omega )}^2 ) \\
 &\quad \quad \times \,(\;\left\| {\nabla  {\overline {\omega }}_1^k }
\right\|_{L^2(\Omega )}^2 + \left\| {\nabla  {\overline
{\omega }}_2^k } \right\|_{L^2(\Omega )}^2 + \left\| {\nabla
 {\overline {\omega }}_3^k }
\right\|_{L^2(\Omega )}^2 ) \\
 \end{split}
\end{equation*}
\\

A convolution inequality in [1] is applied to get
\\
\begin{equation*}
\begin{split}
 &\left\| {\overline {\overline u}_i^k } \right\|_{L^2(\Omega )}^2 \le \left\|
{J_\varepsilon * \left| {\overline u_i^k } \right|(x)}
\right\|_{L^2({\mathbb R}^3)}^2 \le \left\| {J_\varepsilon }
\right\|_{L^1({\mathbb R}^3)}^2 \left\|
{\overline u_i^k } \right\|_{L^2(\Omega )}^2 \\
 &\quad \quad \quad \;\; = \left\| {\overline u_i^k } \right\|_{L^2(\Omega
)}^2 \le \mathop {\sup }\limits_{(t_{k - 1} ,t_k )} \left\| {\tilde
u_i } \right\|_{L^2(\Omega )}^2 \\
 \end{split}
\end{equation*}
Thus we have
\begin{equation*}
\begin{split}
 &\int_\Omega {(\tilde {\omega
}_1^2 + \tilde {\omega }_2^2 + \tilde {\omega }_3^2 )} \;\; +
\int_{t_{k - 1} }^{t } {(\;\left\| {\nabla \tilde {\omega }_1 }
\right\|_{L^2(\Omega )}^2 + \left\| {\nabla \tilde {\omega }_2 }
\right\|_{L^2(\Omega )}^2 +
\left\| {\nabla \tilde {\omega }_3 } \right\|_{L^2(\Omega )}^2 )} \\
 &\qquad \le \int_\Omega {(\tilde {\omega }_1^{k - 1^2} + \tilde {\omega
}_2^{k - 1^2} + \tilde {\omega }_3^{k - 1^2} )} \;\; + \\
 &\qquad \quad + \,C_2\,\left( \;{\mathop {\sup }\limits_{(t_{k - 1} ,t_k )} \{\,\left\| {u_1
} \right\|_{L^2(\Omega )}^2 + \left\| {u_2 } \right\|_{L^2(\Omega
)}^2 +
\left\| {u_3 } \right\|_{L^2(\Omega )}^2 \} \;\;+ } \right. \\
 &\quad \quad \quad \quad \quad \quad \quad  + \left.
{\frac{1}{\Delta t_k }\int_{t_{k - 1} }^{t_k } {\{\,\left\| {\nabla
u_1 } \right\|_{L^2(\Omega )}^2 + \left\| {\nabla u_2 }
\right\|_{L^2(\Omega )}^2
+ \left\| {\nabla u_3 } \right\|_{L^2(\Omega )}^2 \}} } \right) \\
\end{split}
\end{equation*}
\begin{equation*}
\begin{split}
 &\quad \qquad \times \int_{t_{k - 1} }^{t } { \int_\Omega {(\tilde {\omega
}_1^2 + \tilde {\omega }_2^2 + \tilde {\omega }_3^2 )} } \;\; + \\
 &\quad \qquad  + \,C_3\,\left( {\Delta t_k \mathop {\sup
}\limits_{(t_{k - 1} ,t_k )} \{\,\left\| {u_1 } \right\|_{L^2(\Omega
)}^2 + \left\| {u_2 } \right\|_{L^2(\Omega )}^2 + \left\| {u_3 }
\right\|_{L^2(\Omega )}^2 \} \;\;+ } \right. \\
 &\quad \quad \quad \quad \quad \quad \quad  + \left.
{\int_{t_{k - 1} }^{t_k } {\{\,\left\| {\nabla u_1 }
\right\|_{L^2(\Omega )}^2 + \left\| {\nabla u_2 }
\right\|_{L^2(\Omega )}^2 + \left\| {\nabla u_3
} \right\|_{L^2(\Omega )}^2 \}} } \right) {M_{k-1}} \\
 \end{split}
\end{equation*}
where $C_2, C_3 >0$ are constants indepentent of $k$. Set
\begin{equation*}
\begin{split}
 &K_k^ * = \Delta t_k \mathop {\sup }\limits_{(t_{k - 1} ,t_k )} \{\,\left\|
{u_1 } \right\|_{L^2(\Omega )}^2 + \left\| {u_2 }
\right\|_{L^2(\Omega )}^2
+ \left\| {u_3 } \right\|_{L^2(\Omega )}^2 \} + \\
 &\quad \quad \quad \quad \quad + \int_{t_{k - 1} }^{t_k } {\{\,\left\|
{\nabla u_1 } \right\|_{L^2(\Omega )}^2 + \left\| {\nabla u_2 }
\right\|_{L^2(\Omega )}^2 + \left\| {\nabla u_3 }
\right\|_{L^2(\Omega )}^2
\}} \\
 \end{split}
\end{equation*}
and
\[
 f_k (t) = \,\mathop {\sup }\limits_{(t_{k - 1} ,t)}
\int_\Omega {(\tilde {\omega }_1^2 + \tilde {\omega }_2^2 + \tilde
{\omega
}_3^2 )}
\]
Then we arrive at
\[
f_k (t) \le M_{k-1} \;+ C_2 \frac{1}{\Delta t_k } K_k^ *
\;\int_{t_{k-1}}^{t} {f_k (t)}  + \; C_3 K_k^ * {M_{k-1}}
\]
\\
By using Gronwall inequality it follows that
\[
M_k \le \left( 1 + C_3 K_k^ * \right)\; \exp \left( C_2 K_k^ *
\right) {M_{k-1}}
\]
\\

Note that
\begin{equation*}
\begin{split}
 &\sum\limits_{k = 1}^N {K_k^ * } \le T\,\mathop {\sup }\limits_{t \in (0,T)}
\int_\Omega {(u_1^2 + u_2^2 + u_3^2 )} + \int_0^T {(\,\left\|
{\nabla u_1 } \right\|_{L^2(\Omega )}^2 + \left\| {\nabla u_2 }
\right\|_{L^2(\Omega )}^2
+ \left\| {\nabla u_3 } \right\|_{L^2(\Omega )}^2 )} \\
 &\quad \qquad  \le (T+1)\,K_0 < + \infty \\
 \end{split}
\end{equation*}
\\
Hence, combining above two cases, we obtain
\\
\begin{equation*}
\begin{split}
 & M_1 \le \left( 1 + C_3 K_1^ * \right) \;\exp
\left( C_2 K_1^ * \right) \,M_0 \\
 & M_2 \le \left( 1 + C_3 K_1^ * \right)\,\left( 1 + C_3 K_2^ * \right) \;\exp
\left( C_2 \sum\limits_{k = 1}^2 {K_k^ * } \right) \,M_0 \\
 & \; \cdots \; \cdots \cdots \cdots \cdots \cdots \cdots \\
 & M_N \le \prod\limits_{k = 1}^N {(1 + C_3 K_k^\ast )} \;\exp
\left( C_2 \sum\limits_{k = 1}^N {K_k^ * } \right) \,M_0 \\
 \end{split}
\end{equation*}
Note that
\\
\begin{equation*}
\begin{split}
 &\quad \prod\limits_{k = 1}^N {(1 + C_3 K_k^\ast )} = \exp \left( {\ln
\prod\limits_{k = 1}^N {(1 + C_3 K_k^\ast )} } \right)   \\
 &= \exp \left(
{\sum\limits_{k = 1}^N {\ln (1 + C_3 K_k^\ast )} } \right) \le \exp \left(
{C_3 \sum\limits_{k = 1}^N {K_k^\ast } } \right) = \exp (C_3 (T + 1)K_0 )   \\
 \end{split}
\end{equation*}
These mean that
\[
 M_k \le  M_0 \;\exp \left( (C_2 + C_3)\, (T+1)K_0 \right) \qquad   k=1,\cdots ,N
\]

Finally we get
\begin{equation*}
\begin{split}
&\mathop {\sup }\limits_{t \in (0,T)} \int_{\mathbb R^3} {(\tilde {\omega
}_1^2 + \tilde {\omega }_2^2 + \tilde {\omega }_3^2 )} \le
\, \mathop {\max }\limits_k \{M_k \} \\
&\qquad \qquad  \le  M_0 \;\exp \left( (C_2 + C_3) \,(T+1)K_0 \right)   \\
 \end{split}
\end{equation*}

This conclusion is also true for the weak solution of problem (5)
by means of the result of section 3 and the lower limit of Galerkin
sequence according to the page 196 of [4].
\\
\\
\\

\textbf{3. Existence}

In this section we have to consider the existence of solutions of
the auxiliary problems. We just need to consider the following
system on $(0,\delta )$ :
\begin{equation}
\label{eq7}
\begin{split}
 &\partial _t \tilde \omega _1 \,+ \overline{\overline u}_1 \partial _{x_1 } \overline {\omega }_1 + \overline{\overline u}_2
\partial _{x_2 } \overline {\omega }_1 + \overline{\overline u}_3 \partial _{x_3 } \overline {\omega }_1
-\overline {\omega }_1 \partial _{x_1 } \overline{\overline u}_1
-\overline {\omega }_2
\partial _{x_2 } \overline{\overline u}_1 -\overline {\omega }_3 \partial _{x_3 } \overline{\overline u}_1
+\partial _{x_1 } q = \Delta \tilde \omega
_1 \\
 &\partial _t \tilde \omega _2 + \overline{\overline u}_1 \partial _{x_1 } \overline {\omega }_2 + \overline{\overline u}_2 \partial
_{x_2 } \overline {\omega }_2 + \overline{\overline u}_3 \partial
_{x_3 } \overline {\omega }_2 -\overline {\omega }_1 \partial _{x_1
} \overline{\overline u}_2 -\overline {\omega }_2
\partial _{x_2 } \overline{\overline u}_2 -\overline {\omega }_3 \partial _{x_3 } \overline{\overline u}_2
+\partial _{x_2 } q = \Delta \tilde \omega _2
\\
 &\partial _t \tilde \omega _3 + \overline{\overline u}_1 \partial _{x_1 } \overline {\omega }_3 + \overline{\overline u}_2 \partial
_{x_2 } \overline {\omega }_3 + \overline{\overline u}_3 \partial
_{x_3 } \overline {\omega }_3 -\overline {\omega }_1 \partial _{x_1
} \overline{\overline u}_3 -\overline {\omega }_2
\partial _{x_2 } \overline{\overline u}_3 -\overline {\omega }_3 \partial _{x_3 } \overline{\overline u}_3
+\partial _{x_3 } q = \Delta \tilde \omega _3
\\
 \end{split}
\end{equation}
with initial value $\tilde \omega _i (x,0)=\omega _{i0} (i=1,2,3)$ and
\[
\overline {\omega }_i (x)=\frac{1}{\delta }\int_0^\delta {\tilde \omega _i
(x,t)dt}
\]
and
\[
\overline u_i (x)=\frac{1}{\delta }\int_0^\delta u _i (x,t)dt,\quad
\overline{\overline u}_i (x) = J_\varepsilon * \overline u _i (x),
\quad i=1,2,3
\]
as well as the incompressible conditions:
\[
\partial _{x_1 } \tilde \omega _1 + \partial _{x_2 } \tilde \omega _2 + \partial _{x_3 }
\tilde \omega _3 = 0\quad \Rightarrow \quad \partial _{x_1 } \overline
{\omega }_1 +\partial _{x_2 } \overline {\omega }_2 +\partial _{x_3
} \overline {\omega }_3 =0
\]
\[
\partial _{x_1 } u_1 +\partial _{x_2 } u_2 +\partial _{x_3 }
u_3 =0\quad \Rightarrow \quad \partial _{x_1 } \overline{\overline
u}_1 +\partial _{x_2 } \overline{\overline u}_2 +\partial _{x_3 }
\overline{\overline u}_3 =0
\]

(i) The Galerkin procedure is applied. For each $m$ and $i=1,2,3$ we
define an approximate solution $(\tilde \omega _{1m} ,\;\tilde \omega _{2m}
,\;\tilde \omega _{3m} )$ as follows:
\[
\tilde \omega _{im} =\sum\limits_{j=1}^m {g_{ij}(t) w_{ij} }
\]
where $\{w_{i1} ,\;\cdots ,\;w_{im} ,\cdots \}$ is the basis of $W$,
and $W=$ the closure of $\mathscr{V}$ in the Sobolev space $
W^{2,4}(\Omega) $, which is separable and is dense in $V$. Thus
\begin{equation}
\label{eq8}  {(\partial _t \tilde \omega _{im} ,\;w_{il} )} + {(\nabla
\tilde \omega _{im} ,\;\nabla w_{il} )} + {((\overline{\overline u}\cdot
\nabla )\overline {\omega }_{im} ,\;w_{il} )} - {((\overline {\omega
}_m \cdot \nabla ) \overline{\overline u}_i ,\;w_{il} )} =0
\end{equation}
\begin{equation*}
\begin{split}
 &t\in (0,\delta ),\quad \quad l=1,\cdots ,m \\
 &\\
 &\tilde \omega _{im} (0) = \omega _{i0}^m \\
 \end{split}
\end{equation*}
where $\omega _{i0}^m $ is the orthogonal projection in $H$ of
$\omega _{i0} $ onto the space spanned by $w_{i1} ,\;\cdots
\;w_{im}$. Therefore,
\begin{equation*}
\begin{split}
 &\sum\limits_{j=1}^m { {(w_{ij} ,\;w_{il} ){g}'_{ij} (t)}
} +\sum\limits_{j=1}^m { {(\nabla w_{ij} ,\;\nabla w_{il}
)g_{ij} (t)} } + \\
 &\quad +\sum\limits_{j=1}^m { {\{(( \overline{\overline u}(t)\cdot \nabla
)w_{ij} ,\;w_{il} )-((w_j \cdot \nabla )
w_{il},\;\overline{\overline u}_i (t) )\}}
}\; \overline {g}_{ij} (t)\, =0 \\
 \end{split}
\end{equation*}
where $\overline {g}_{ij} (t)= \frac{1}{\delta } \int_0^\delta {
g_{ij} (t)dt}$, and $ {u}_i \in L^{\infty}(0,T; H)$ from section 1 which are determined by equations (1).
Inverting the nonsigular matrix with elements $ {(w_{ij} ,\;w_{il}
)} ,\;\;1\le j,l\le m$, we can write above system in the following
form
\begin{equation}
\label{eq9} {g}'_{ij} (t)+\sum\limits_{l=1}^m { {\alpha _{ijl}
g_{il} (t)} } +\sum\limits_{l=1}^m { {\beta _{ijl}\; \overline
{g}_{il} (t)} }=0
\end{equation}
where $\alpha _{ijl} ,\;\,\beta_{ijl}$ are constants.

The initial conditions are equivalent to
\[
g_{ij}
(0)=g^0_{ij}=\mbox{the}\;j^{th}\;\mbox{component}\;\mbox{of}\;\omega
_{i0}^m
\]

We construct a sequence $\{g_{ij}^k \}$ by using a successive
approximation:
\begin{equation*}
\begin{split}
 &{g_{ij}^1} ^\prime =-\sum\limits_{l=1}^m {\alpha _{ijl} g_{il}^0 }
-\sum\limits_{l=1}^m {\beta _{ijl} \overline {g}_{il}^0 } \quad
\Rightarrow \quad g_{ij}^1 =g_{ij}^0 -\int_0^t {\left(
{\sum\limits_{l=1}^m {\alpha _{ijl} g_{il}^0 } +\sum\limits_{l=1}^m
{\beta _{ijl} \overline {g}_{il}^0 } } \right)}
\\
 &{g_{ij}^2} ^\prime =-\sum\limits_{l=1}^m {\alpha _{ijl} g_{il}^1 }
-\sum\limits_{l=1}^m {\beta _{ijl} \overline {g}_{il}^1 } \quad
\Rightarrow \quad g_{ij}^2 =g_{ij}^0 -\int_0^t {\left(
{\sum\limits_{l=1}^m {\alpha _{ijl} g_{il}^1 } +\sum\limits_{l=1}^m
{\beta _{ijl} \overline {g}_{il}^1 } } \right)}
\\
 &\quad \quad \quad \quad \cdots \cdots \cdots \cdots \\
 &{g_{ij}^k} ^\prime =-\sum\limits_{l=1}^m {\alpha _{ijl} g_{il}^{k-1} }
-\sum\limits_{l=1}^m {\beta _{ijl} \overline {g}_{il}^{k-1} } \quad
\Rightarrow \quad g_{ij}^k =g_{ij}^0 -\int_0^t {\left(
{\sum\limits_{l=1}^m {\alpha _{ijl} g_{il}^{k-1} }
+\sum\limits_{l=1}^m {\beta _{ijl} \overline {g}_{il}^{k-1}
} } \right)} \\
 \end{split}
\end{equation*}
so that
\[
\qquad \left| {g_{ij}^k (t)-g_{ij}^{k-1} (t)} \right|\le \int_0^t
{\left( {\sum\limits_{l=1}^m {\left| {\alpha _{ijl} } \right|\left|
{g_{il}^{k-1} (t)-g_{il}^{k-2} (t)} \right|} +\sum\limits_{l=1}^m
{\left| {\beta _{ijl} } \right|\left| {\overline {g}_{il}^{k-1}
(t)-\overline {g}_{il}^{k-2} (t)} \right|} } \right)}
\]
it follows that
\[
\qquad \mathop {\max }\limits_{i,j} \mathop {\sup }\limits_t \left|
{g_{ij}^k (t)-g_{ij}^{k-1} (t)} \right|\le \mathop {\max
}\limits_{i,j} \sum\limits_{l=1}^m {\left( {\left| {\alpha _{ijl} }
\right|+\left| {\beta _{ijl} } \right|} \right)} \;\cdot\;
t\;\cdot\; \mathop {\max }\limits_{i,j} \mathop {\sup }\limits_t
\left| {g_{ij}^{k-1} (t)-g_{ij}^{k-2} (t)} \right|
\]

Taking $ \delta :\,=\frac{1}{\mathop {\max }\limits_{i,j}
\sum\limits_{l=1}^m {\left( {\left| {\alpha _{ijl} } \right|+2\left|
{\beta _{ijl} } \right|} \right)} } $, as $t\le \delta $, then
choosing $\delta ^\ast $:
\[
0<\delta ^\ast =\frac{\mathop {\max }\limits_{i,j}
\sum\limits_{l=1}^m {\left( {\left| {\alpha _{ijl} } \right|+\left|
{\beta _{ijl} } \right|} \right)} }{\mathop {\max }\limits_{i,j}
\sum\limits_{l=1}^m {\left( {\left| {\alpha _{ijl} } \right|+2\left|
{\beta _{ijl} } \right|} \right)} }<1
\]
we have
\[
\qquad \mathop {\max }\limits_{i,j} \left\| {g_{ij}^k -g_{ij}^{k-1}
} \right\|_\infty \le \delta ^\ast \cdot \mathop {\max
}\limits_{i,j} \left\| {g_{ij}^{k-1} -g_{ij}^{k-2} } \right\|_\infty
\le (\delta ^\ast )^{k-1}\cdot \mathop {\max }\limits_{i,j} \left\|
{g_{ij}^1 -g_{ij}^0 } \right\|_\infty
\]

For any $n,k$ (we can set $n>k$ without loss of generality), we get
\begin{equation*}
\begin{split}
 &\qquad \mathop {\max }\limits_{i,j} \left\| {g_{ij}^n -g_{ij}^k } \right\|_\infty
\le \mathop {\max }\limits_{i,j} \left\| {g_{ij}^n -g_{ij}^{n-1} }
\right\|_\infty +\cdots +\mathop {\max }\limits_{i,j} \left\|
{g_{ij}^{k+1}
-g_{ij}^k } \right\|_\infty \\
 &\qquad \le \left( {(\delta ^\ast )^{n-1}+\cdots +(\delta ^\ast )^k} \right)\cdot
\mathop {\max }\limits_{i,j} \left\| {g_{ij}^1 -g_{ij}^0 }
\right\|_\infty =(\delta ^\ast )^k\frac{1-(\delta ^\ast
)^{n-k}}{1-\delta ^\ast }\mathop
{\max }\limits_{i,j} \left\| {g_{ij}^1 -g_{ij}^0 } \right\|_\infty \\
 &\qquad \to 0\quad (k\to \infty ) \\
 \end{split}
\end{equation*}
Thus, for every $i=1,2,3;\;\;j=1,\cdots ,m$, $\{g_{ij}^k \}$ is a
Cauchy sequence in $L^\infty (0,\delta )$. Since $L^\infty (0,\delta
)$ is complete, then there exists a function $g_{ij}^\ast \in
L^\infty (0,\delta )$ such that $\left\| {g_{ij}^k -g_{ij}^\ast }
\right\|_\infty \to 0$ as $k\to \infty $.
\\

From
\[
g_{ij}^k (t)=g_{ij}^0 -\int_0^t {\left( {\sum\limits_{l=1}^m {\alpha
_{ijl} g_{il}^{k-1} (t)} +\sum\limits_{l=1}^m {\beta _{ijl}
\overline {g}_{il}^{k-1} (t)} } \right)}
\]
let $k\to \infty $, it follows that
\[
g_{ij}^\ast (t)=g_{ij}^0 -\int_0^t {\left( {\sum\limits_{l=1}^m
{\alpha _{ijl} g_{il}^\ast (t)} +\sum\limits_{l=1}^m {\beta _{ijl}
\overline {g}_{il}^\ast (t)} } \right)}
\]
i.e., $g_{ij}^\ast $ is a solution of the system (9) on $(0,\delta
)$ for which $g_{ij}^\ast (0)=g_{ij}^0 $, $i=1,2,3; \;\, j=1,\cdots
,m$.
\\

 (ii)
\[
\sum\limits_{i=1}^3 {(\partial _t \tilde \omega _{im} ,\;\tilde \omega _{im} )}
+\sum\limits_{i=1}^3 {(\nabla \tilde \omega _{im} ,\;\nabla \tilde \omega _{im} )}
+\sum\limits_{i=1}^3 {(( \overline{\overline u} \cdot \nabla
)\overline {\omega }_{im} ,\;\tilde \omega _{im} )} - \sum\limits_{i=1}^3
{((\overline {\omega }_m \cdot \nabla ) \overline{\overline u}_i
,\;\tilde \omega _{im} )} = 0
\]
Then we write
\[
\frac{1}{2}\frac{d}{dt}\left( {\sum\limits_{i=1}^3 {\left\| {\tilde \omega
_{im} } \right\|_{L^2(\Omega )}^2 } } \right)+\sum\limits_{i=1}^3
{\left\| {\nabla \tilde \omega _{im} } \right\|_{L^2(\Omega )}^2 }
-\sum\limits_{i=1}^3 {(( \overline{\overline u}\cdot \nabla )\tilde \omega
_{im} ,\;\overline {\omega }_{im} )} +\sum\limits_{i=1}^3
{((\overline {\omega }_m \cdot \nabla )\tilde \omega _{im} ,\;
\overline{\overline u}_i )} = 0
\]

Similar to those in the section 2, and $\eta $ is chosen to be small
enough, we have
\[
\sum\limits_{i=1}^3 {\left\| {\tilde \omega _{im} } \right\|_{L^2(\Omega
)}^2 } + \int_0^\eta {\left( {\sum\limits_{i=1}^3 {\left\| {\nabla
\tilde \omega _{im} } \right\|_{L^2(\Omega )}^2 } } \right)} \le 2 \left(
{\sum\limits_{i=1}^3 {\left\| {\omega _{i0}^m } \right\|_{L^2(\Omega
)}^2 } } \right)
\]
as $ 1-4\,K_0\,\eta/\mu_\varepsilon  \ge 1/2 $ . Hence
\begin{equation}
\label{eq10} \mathop {\sup }\limits_{t\in (0,\eta )} \left(
{\sum\limits_{i=1}^3 {\left\| {\tilde \omega _{im} } \right\|_{L^2(\Omega
)}^2 } } \right)\le 2 \left( {\sum\limits_{i=1}^3 {\left\| {\omega
_{i0} } \right\|_{L^2(\Omega )}^2 } } \right)
\end{equation}
and
\begin{equation}
\label{eq11} \sum\limits_{i=1}^3 {\left\| {\tilde \omega _{im} (\eta )}
\right\|_{L^2(\Omega )}^2 } + \int_0^\eta {\sum\limits_{i=1}^3
{\left\| {\nabla \tilde \omega _{im} } \right\|_{L^2(\Omega )}^2 } } \le 2
\,\left( {\sum\limits_{i=1}^3 {\left\| {\omega _{i0} }
\right\|_{L^2(\Omega )}^2 } } \right)
\end{equation}

The inequalities (10) and (11) are valid for any fixed $\delta \le
\eta$.
\\

(iii) Let $ {\mathop \omega \limits^\circ }_m $ denote the function from ${\mathbb
R}$ into $V$, which is equal to $\tilde \omega _m $ on $(0,\delta )$ and to
0 on the complement of this interval. The Fourier transform of
$ {\mathop \omega \limits^\circ }_m $ is denoted by $\hat {\omega }_m $. We want to
show that
\[
\int_{-\infty }^{+\infty } {\left| \tau \right|^{2\gamma }\left(
{\sum\limits_{i=1}^3 {\left\| {\hat {\omega }_{im} (\tau )}
\right\|_{L^2(\Omega )}^2 } } \right)} \,d\tau <+\infty
\]
for some $\gamma >0$. Along with (\ref{eq11}) this will imply that
\\

$ {\mathop \omega \limits^\circ }_m $ belongs to a bounded set of $H^\gamma
({\mathbb R},V,H)\\$
\\
and will enable us to apply the result of compactness.

We observe that (\ref{eq8}) can be written as
\[
\frac{d}{dt}\left( {\sum\limits_{i=1}^3 {( {\mathop \omega \limits^\circ }_{im}
,\;w_{ij} )} } \right)=\sum\limits_{i=1}^3 {( {\mathop f \limits^\circ }_{im}
,\;w_{ij} )} + \sum\limits_{i=1}^3 {(\omega _{i0}^m ,\;w_{ij} )\,}
\eta _0 -\sum\limits_{i=1}^3 {(\tilde \omega _{im} (\delta ),\;w_{ij} )\,}
\eta _\delta
\]
where $\eta _0 ,\;\eta _\delta $ are Dirac distributions at 0 and
$\delta $, and
\[
f_{im} = -\Delta \tilde \omega _{im} +( \overline{\overline u} \cdot \nabla
)\overline {\omega }_{im} -(\overline {\omega }_m \cdot \nabla )
\overline{\overline u}_i\;\;\;
\]
\[
{\mathop f \limits^\circ }_{im} = f_{im} \;\mbox{on}\;(0,\delta ),\quad 0\;
\mbox{outside this interval}\]

By the Fourier transform,
\[
2\mbox{i}\pi \tau \sum\limits_{i=1}^3 {(\hat {\omega }_{im}
,\;w_{ij} )} =\sum\limits_{i=1}^3 {(\hat {f}_{im} ,\;w_{ij} )}
+\sum\limits_{i=1}^3 {(\omega _{i0}^m ,\;w_{ij} )}
-\sum\limits_{i=1}^3 {(\tilde \omega _{im} (\delta ),\;w_{ij} )\,} \exp
(-2\mbox{i}\pi \delta \tau )
\]
where $\hat {\omega }_{im} $ and $\hat {f}_{im} $ denoting the
Fourier transforms of ${\mathop \omega \limits^\circ }_{im} $ and ${\mathop f \limits^\circ }_{im}
$ respectively.
\\

We multiply above equality by $\hat {g}_{ij} (\tau )=$Fourier
transform of ${\mathop g \limits^\circ }_{ij} $ and add the resulting equation for
$j=1,\cdots ,m$, we get
\begin{equation*}
\begin{split}
&2\mbox{i}\pi \tau \sum\limits_{i=1}^3 {\left\| {\hat {\omega }_{im}
(\tau )} \right\|_{L^2(\Omega )}^2 } =\sum\limits_{i=1}^3 {(\hat
{f}_{im} (\tau ),\;\hat {\omega }_{im} (\tau ))}\\
 &\quad \quad
+\sum\limits_{i=1}^3 {(\omega _{i0}^m ,\;\hat {\omega }_{im} (\tau
))} -\sum\limits_{i=1}^3 {(\tilde \omega _{im} (\delta ),\;\hat {\omega
}_{im} (\tau ))\,} \exp (-2\mbox{i}\pi \delta \tau )
\end{split}
\end{equation*}

For some $\varphi _i \in V$,
\begin{equation*}
\begin{split}
 &\int_0^\delta {\sum\limits_{i=1}^3 {(f_{im} ,\;\varphi _i )} }
=\int_0^\delta {\sum\limits_{i=1}^3 {(-\Delta \tilde \omega _{im}
,\;\varphi _i )} } +\int_0^\delta {\sum\limits_{i=1}^3 {((
\overline{\overline u} \cdot \nabla )\overline {\omega }_{im}
,\;\varphi _i )} } -\int_0^\delta {\sum\limits_{i=1}^3 {((\overline
{\omega }_m
\cdot \nabla ) \overline{\overline u}_i ,\;\varphi _i )} } \\
 &\qquad =\int_0^\delta {\sum\limits_{i=1}^3 {(\nabla \tilde \omega _{im} ,\;\nabla \varphi
_i )} } -\int_0^\delta {\sum\limits_{i=1}^3 {(( \overline{\overline
u} \cdot \nabla )\varphi _i ,\;\overline {\omega }_{im} )} }
+\int_0^\delta {\sum\limits_{i=1}^3 {((\overline
{\omega }_m \cdot \nabla )\varphi _i ,\; \overline{\overline u}_i )} } \\
 &\qquad \le \int_0^\delta {\sum\limits_{i=1}^3 {\left\| {\nabla \tilde \omega _{im} }
\right\|_{L^2(\Omega )} \left\| {\nabla \varphi _i }
\right\|_{L^2(\Omega )}
} } \;\; + \\
 &\qquad \qquad + 2 \int_0^\delta \left( {\sum\limits_{i=1}^3 {\left\| { \overline{\overline u}_i } \right\|_{L^{\infty}(\Omega
)}^2 } } \right)^{1/2}\left( {\sum\limits_{i=1}^3 {\left\|
{\overline {\omega }_{im} } \right\|_{L^2(\Omega )}^2 } }
\right)^{1/2}\left( {\sum\limits_{i=1}^3 {\left\| {\nabla \varphi _i
} \right\|_{L^2(\Omega
)}^2 } } \right)^{1/2} \\
 \end{split}
\end{equation*}
\begin{equation*}
\begin{split}
 &\qquad \le \int_0^\delta {\left( {\sum\limits_{i=1}^3 {\left\| {\nabla \tilde \omega
_{im} } \right\|_{L^2(\Omega )}^2 } } \right)^{1/2}} \left(
{\sum\limits_{i=1}^3 {\left\| {\nabla \varphi _i }
\right\|_{L^2(\Omega )}^2
} } \right)^{1/2} \; + \\
 &\qquad \qquad + 2\,\delta \left( {\sum\limits_{i=1}^3 {\left\| { \overline u_i } \right\|_{L^2(\Omega
)}^2 } } \right)^{1/2}\left( {\sum\limits_{i=1}^3 {\left\|
{\overline {\omega }_{im} } \right\|_{L^2(\Omega )}^2 } }
\right)^{1/2}\left( {\sum\limits_{i=1}^3 {\left\| {\nabla \varphi _i
} \right\|_{L^2(\Omega
)}^2 } } \right)^{1/2} \\
 &\qquad \le \int_0^\delta {\left( {\sum\limits_{i=1}^3 {\left\| {\nabla \tilde \omega
_{im} } \right\|_{L^2(\Omega )}^2 } } \right)^{1/2}} \left\| {\nabla
\varphi
} \right\|_V \;\; + \\
 &\qquad \qquad + 2 \delta \left( { \mathop {\;\sup
}\limits_{(0,\delta )} \sum\limits_{i=1}^3 {\left\| {u_i }
\right\|_{L^2(\Omega )}^2 } } \right)^{1/2}\; \left( { \mathop
{\;\sup }\limits_{(0,\delta )} \sum\limits_{i=1}^3 {\left\| {\tilde \omega
_{im} } \right\|_{L^2(\Omega )}^2 } } \right)^{1/2}\left\| {\nabla
\varphi }
\right\|_V \\
 \end{split}
\end{equation*}
this remains bounded according to (\ref{eq2}), (\ref{eq3}), and
(\ref{eq10}), (\ref{eq11}). Therefore
\[
\int_0^\delta {\left\| {f_{im} (t)} \right\|_V dt} =\int_0^\delta
{\;\mathop {\sup }\limits_{\left\| \varphi \right\|_V =1}
\;\sum\limits_{i=1}^3 {(f_{im} ,\;\varphi _i )} } <+\infty
\]
it follows that
\[
\mathop {\sup }\limits_{\tau \in {\mathbb R}} \left\| {\hat {f}_{im}
(\tau )} \right\|_V <+\infty ,\quad \;\forall m
\]

Due to (\ref{eq10}), we have
\[
\left\| {\omega _{im} (0)} \right\|_{L^2(\Omega )} <+\infty ,\quad
\quad \left\| {\tilde \omega _{im} (\delta )} \right\|_{L^2(\Omega )}
<+\infty
\]
then by Poincare inequality,
\begin{equation}
\label{eq12}
\begin{split}
 \left| \tau \right|\sum\limits_{i=1}^3 {\left\| {\hat
{\omega }_{im} (\tau )} \right\|_{L^2(\Omega )}^2 } & \le c_1
\sum\limits_{i=1}^3 {\left\| {\hat {f}_{im} (\tau )} \right\|_V
\;\left\| {\hat {\omega }_{im} (\tau )} \right\|_V } + c_2
\sum\limits_{i=1}^3 {\left\| {\hat {\omega }_{im} (\tau )}
\right\|_{L^2(\Omega )} }   \\
 &\le c_3 \left( \sum\limits_{i=1}^3 {\left\| { \hat {\omega }_{im}
(\tau )} \right\|_{L^2(\Omega )}} + \sum\limits_{i=1}^3 {\left\|
{\nabla \hat {\omega }_{im} (\tau )} \right\|_{L^2(\Omega )}}
\right)   \\
\end{split}
\end{equation}

For $\gamma $ fixed, $\gamma <1/4$, we observe that
\[
\left| \tau \right|^{2\gamma }\le c_4 (\gamma )\frac{1+\left| \tau
\right|}{1+\left| \tau \right|^{1-2\gamma }},\quad \quad \forall
\tau \in {\mathbb R}
\]
Thus by (\ref{eq12}),
\begin{equation*}
\begin{split}
 &\int_{-\infty }^{+\infty } {\left| \tau \right|^{2\gamma }\left(
{\sum\limits_{i=1}^3 {\left\| {\hat {\omega }_{im} (\tau )}
\right\|_{L^2(\Omega )}^2 } } \right)} \,d\tau \le c_4 (\gamma
)\int_{-\infty }^{+\infty } {\frac{1+\left| \tau \right|}{1+\left|
\tau \right|^{1-2\gamma }}\left( {\sum\limits_{i=1}^3 {\left\| {\hat
{\omega }_{im} (\tau )} \right\|_{L^2(\Omega )}^2 } } \right)}
\,d\tau   \\
 &\le c_5 \int_{-\infty }^{+\infty } {\frac{1}{1+\left| \tau
\right|^{1-2\gamma }}\sum\limits_{i=1}^3 {\left\| { \hat
{\omega }_{im} (\tau )} \right\|_{L^2(\Omega )} } } d\tau \\
 & + \;c_6 \int_{-\infty }^{+\infty } {\frac{1}{1+\left| \tau
\right|^{1-2\gamma }}\sum\limits_{i=1}^3 {\left\| {\nabla \hat
{\omega }_{im} (\tau )} \right\|_{L^2(\Omega )} } } d\tau +
 \;c_7 \int_{-\infty }^{+\infty } {\sum\limits_{i=1}^3 {\left\| {
\hat {\omega }_{im} (\tau )} \right\|_{L^2(\Omega )}^2 } } d\tau  \\
\end{split}
\end{equation*}
Because of the Parseval equality,
\begin{equation*}
\begin{split}
 & \int_{-\infty }^{+\infty } {\sum\limits_{i=1}^3 {\left\| {\hat {\omega
}_{im} (\tau )} \right\|_{L^2(\Omega )}^2 d\tau } } = \int_0^\delta
{\sum\limits_{i=1}^3 {\left\| {\tilde \omega _{im} (t)}
\right\|_{L^2(\Omega )}^2
dt} } \\
 &\qquad \qquad \qquad \quad \quad \quad \quad \quad \;\;\, \le \,\delta \;\mathop {\sup
}\limits_{(0,\delta )} \;\sum\limits_{i=1}^3 {\left\| {\tilde \omega _{im}
}
\right\|_{L^2(\Omega )}^2 } <+\infty \\
 & \int_{-\infty }^{+\infty } {\sum\limits_{i=1}^3 {\left\| {\nabla
\hat {\omega }_{im} (\tau )} \right\|_{L^2(\Omega )}^2 } } \,d\tau
=\int_0^\delta {\sum\limits_{i=1}^3 {\left\| {\nabla \tilde \omega _{im}
(t)} \right\|_{L^2(\Omega )}^2 } } \,dt<+\infty  \\
\end{split}
\end{equation*}
as $m\to \infty $. By Cauchy-Schwarz inequality and the Parseval
equality,
\begin{equation*}
\begin{split}
 &\int_{-\infty }^{+\infty } {\frac{1}{1+\left| \tau \right|^{1-2\gamma
}}\sum\limits_{i=1}^3 {\left\| { \hat {\omega }_{im} (\tau )}
\right\|_{L^2(\Omega )} } } d\tau \\
 &\qquad \le \left( {\int_{-\infty }^{+\infty } {\frac{1}{(1+\left| \tau
\right|^{1-2\gamma })^2}} } \right)^{1/2}\left( {\int_0^\delta
{\sum\limits_{i=1}^3 {\left\| { \tilde \omega _{im} (t)}
\right\|_{L^2(\Omega
)}^2 } } dt} \right)^{1/2}<+\infty \\
 &\int_{-\infty }^{+\infty } {\frac{1}{1+\left| \tau \right|^{1-2\gamma
}}\sum\limits_{i=1}^3 {\left\| {\nabla \hat {\omega }_{im} (\tau )}
\right\|_{L^2(\Omega )} } } d\tau \\
 &\qquad \le \left( {\int_{-\infty }^{+\infty } {\frac{1}{(1+\left| \tau
\right|^{1-2\gamma })^2}} } \right)^{1/2}\left( {\int_0^\delta
{\sum\limits_{i=1}^3 {\left\| {\nabla \tilde \omega _{im} (t)}
\right\|_{L^2(\Omega
)}^2 } } dt} \right)^{1/2}<+\infty \\
 \end{split}
\end{equation*}
as $m\to \infty $ by $\gamma <1/4$ and (\ref{eq11}).
\\

(iv) The estimate (\ref{eq10}) and (\ref{eq11}) enable us to assert
the existence of an element $\tilde \omega ^\ast \in L^2(0,\delta ;V) \cap
L^{\infty}(0,\delta ;H) $ and a subsequence $\tilde \omega _{{m}'} $ such
that
\\

$\tilde \omega _{{m}'} \to \tilde \omega ^\ast $ in $L^2(0,\delta ;V)$ weakly, and
in $L^{\infty}(0,\delta ;H)$ weak-star, as ${m}'\to \infty $
\\

Due to (iii) we also have
\\

$\tilde \omega _{{m}'} \to \tilde \omega ^\ast $ in $L^2(0,\delta ;H)$ strongly as
${m}'\to \infty \\$
\\
This convergence result enable us to pass to the
limit.
\\

Let $\psi _i $ be a continuously differentiable function on
$(0,\delta )$ with $\psi _i (\delta )=0$. We multiply (\ref{eq8}) by
$\psi _i (t)$ then integrate by parts. This leads to the equation
\begin{equation*}
\begin{split}
 &-\int_0^\delta {\sum\limits_{i=1}^3 {(\tilde \omega _{im} (t),\;\partial _t \psi
_i (t)w_{ij} )\,dt} } + \int_0^\delta {\sum\limits_{i=1}^3 {(\nabla
\tilde \omega
_{im} ,\;\psi _i (t)\nabla w_{ij} )\,dt} } \\
 &+\int_0^\delta {\sum\limits_{i=1}^3 {(( \overline{\overline u} \cdot \nabla )\overline {\omega }_{im}
,\;w_{ij} \psi _i (t))} } -\int_0^\delta {\sum\limits_{i=1}^3
{((\overline {\omega }_m \cdot \nabla ) \overline{\overline u}_i
,\;w_{ij} \psi _i (t))} } =\sum\limits_{i=1}^3
{(\omega _{i0}^m ,\;w_{ij} )\psi _i (0)} \\
 \end{split}
\end{equation*}

Since $\tilde \omega _{i{m}'} $ converges to $\tilde \omega _i^\ast $ in
$L^2(0,\delta ;H)$ strongly as ${m}'\to \infty $, then $\overline
{\omega }_{i{m}'} $ also converges strongly to $\overline \omega
_i^\ast $, and
\begin{equation*}
\begin{split}
 &\int_0^\delta {\sum\limits_{i=1}^3 {(\tilde \omega _{i{m}'} ,\;\partial _t
\psi _i (t)w_{ij} )\,dt} } \to \int_0^\delta {\sum\limits_{i=1}^3
{(\tilde \omega _i^\ast ,\;\partial _t \psi _i (t)w_{ij} )\,dt} }\\
 &\int_0^\delta {\sum\limits_{i=1}^3 {(\nabla \tilde \omega _{i{m}'} ,\;\psi _i
(t)\nabla w_{ij} )\,dt} } = -\int_0^\delta {\sum\limits_{i=1}^3
{(\tilde \omega
_{i{m}'} ,\;\psi _i (t)\Delta w_{ij} )\,dt} } \\
 &\quad \;\to -\int_0^\delta {\sum\limits_{i=1}^3 {(\tilde \omega _i^\ast ,\;\psi _i
(t)\Delta w_{ij} )} } =\int_0^\delta {\sum\limits_{i=1}^3 {(\nabla
\tilde \omega
_i^\ast ,\;\psi _i (t)\nabla w_{ij} )\,dt} } \\
\end{split}
\end{equation*}
\begin{equation*}
\begin{split}
 &\int_0^\delta {\sum\limits_{i=1}^3 {(( \overline{\overline u} \cdot \nabla )\overline {\omega }_{i{m}'}
,\;w_{ij} \psi _i (t))} } = -\int_0^\delta {\sum\limits_{i=1}^3 {((
\overline{\overline u}\cdot
\nabla )w_{ij} \psi _i (t),\;\overline {\omega }_{i{m}'} )} } \\
 &\quad \;\to -\int_0^\delta {\sum\limits_{i=1}^3 {(( \overline{\overline u} \cdot \nabla )w_{ij}
\psi _i (t),\;\overline \omega _i^\ast )} } =\int_0^\delta
{\sum\limits_{i=1}^3
{(( \overline{\overline u} \cdot \nabla )\overline \omega _i^\ast ,\;w_{ij} \psi _i (t))} } \\
 &\int_0^\delta {\sum\limits_{i=1}^3 {(( \overline {\omega }_{{m}'} \cdot \nabla ) \overline{\overline u}_i
,\;w_{ij} \psi _i (t))} } \to \int_0^\delta {\sum\limits_{i=1}^3 {((
\overline \omega^\ast \cdot \nabla ) \overline{\overline u}_i,\;w_{ij} \psi _i (t))} } \\
 &\sum\limits_{i=1}^3 {(\omega _{i0}^{{m}'} ,\;w_{ij} )\psi _i (0)}
\to \sum\limits_{i=1}^3 {(\omega _{i0} ,\;w_{ij} )\psi _i (0)}
\\
\end{split}
\end{equation*}
\\
Thus, in the limit we find
\begin{equation}
\label{eq13}
\begin{split}
 &-\int_0^\delta {\sum\limits_{i=1}^3 {(\tilde \omega _i^\ast ,\;\partial _t \psi _i
(t)v_i )\,dt} } +\int_0^\delta {\sum\limits_{i=1}^3 {(\nabla \tilde \omega
_i^\ast
,\;\psi _i (t)\nabla v_i )\,dt} } \\
 &+\int_0^\delta {\sum\limits_{i=1}^3 {(( \overline{\overline u} \cdot \nabla )\overline \omega _i^\ast ,\;v_i
\psi _i (t))} } -\int_0^\delta {\sum\limits_{i=1}^3 {((\overline
\omega ^\ast \cdot \nabla ) \overline{\overline u}_i ,\;v_i \psi _i
(t))} } =\sum\limits_{i=1}^3 {(\omega _{i0}
,\;v_i )\psi _i (0)} \\
 \end{split}
\end{equation}
holds for $v_i =w_{i1} ,\;w_{i2} ,\cdots $; by this equation holds
for $v_i =$any finite linear combination of the $w_{ij} $, and by a
continuity argument above equation is still true for any $v_i \in
V$. Hence we find that $\tilde \omega _i^\ast (i=1,2,3)$ is a Leray-Hopf
weak solution of the system (\ref{eq7}).
\\

Finally it remains to prove that $\tilde \omega _i^\ast $ satisfies the
initial conditions. For this we multiply (\ref{eq7}) by $v_i \psi _i
(t)$, after integrating some terms by parts, we get
\begin{equation*}
\begin{split}
 &-\int_0^\delta {\sum\limits_{i=1}^3 {(\tilde \omega _i^\ast ,\;\partial _t \psi _i
(t)v_i )} } +\int_0^\delta {\sum\limits_{i=1}^3 {(\nabla \tilde \omega
_i^\ast
,\;\psi _i (t)\nabla v_i )\,dt} } \\
 &+\int_0^\delta {\sum\limits_{i=1}^3 {(( \overline{\overline u} \cdot \nabla )\overline \omega _i^\ast ,\;v_i
\psi _i (t))} } -\int_0^\delta {\sum\limits_{i=1}^3 {((\overline
\omega ^\ast \cdot \nabla ) \overline{\overline u}_i ,\;v_i \psi _i
(t))} } =\sum\limits_{i=1}^3 {(\tilde \omega _i^\ast
(0),\;v_i )\psi _i (0)} \\
 \end{split}
\end{equation*}
\\
By comparison with (\ref{eq13}),
\[
\sum\limits_{i=1}^3 {(\tilde \omega _i^\ast (0)-\omega _{i0} ,\;v_i )\psi
_i (0)} =0
\]
Therefore we can choose $\psi _i $ particularly such that
\[
(\tilde \omega _i^\ast (0)-\omega _{i0} ,\;v_i )=0,\quad \quad \forall v_i
\in V
\]
\\
\\
\\

\textbf{4. Convergence}

Now the partition is refined infinitely and $\varepsilon$ becomes
sufficiently small, we will prove that there exists some subsequence
of the solutions of auxiliary problems which converges to a weak
solution of (\ref{eq4}).
\\

Since
\[
\mathop {\sup }\limits_{t\in (0,T)} \int_\Omega {(\tilde {\omega }_1^2
+\tilde {\omega }_2^2 +\tilde {\omega }_3^2 )} \;<+\infty
\]
the family $(\tilde {\omega }_1 ,\tilde {\omega }_2 ,\tilde {\omega
}_3 )$ is uniformly bounded in $L^2 (0,T;H) \cap L^\infty (0,T;H)$,
then we can choose ${\varepsilon}' \to 0 \;
({\varepsilon}'=O(\frac{1}{{k}'}))$ and ${k}'\to \infty $, or
$\Delta t_k ^\prime \to 0$, such that there exists a subsequence
$({\tilde {\omega }}'_1 ,{\tilde {\omega }}'_2 ,{\tilde {\omega
}}'_3 )$ converging weakly in $L^2 (0,T; H)$ and weak-star in
$L^\infty (0,T; H)$ to some element $(\omega _1^\ast ,\omega _2^\ast
,\omega _3^\ast )$. On the other hand, because $\tilde {\omega }_i
(i=1,2,3)$ belong to $L^2(0,T; H)$, we can verify that
\[
\overline {\omega }_i (x,t)=\left\{ {\frac{1}{\Delta t_k
}\int_{t_{k-1} }^{t_k } {\tilde {\omega }_i (x,t)dt} ,\;\;t\in
(t_{k-1} ,t_k )\subset (0,T)} \right\}
\]
also belongs to $L^2(0,T; H)$. In fact,
\begin{equation*}
\begin{split}
 &\int_0^T \int_\Omega {\overline {\omega }_i^2 (x,t)} =\sum\limits_k {\int_{t_{k-1} }^{t_k }
\int_\Omega {\left( {\frac{1}{\Delta t_k }\int_{t_{k-1} }^{t_k }
{\tilde {\omega }_i (x,t)} } \right)} } ^2=\sum\limits_k
{\frac{1}{\Delta t_k^2 }\cdot \Delta t_k \cdot \int_\Omega \left(
{\int_{t_{k-1} }^{t_k } {\tilde {\omega }_i (x,t)} }
\right)} ^2 \\
 &\quad \le \sum\limits_k {\frac{1}{\Delta t_k }\cdot \int_\Omega \left( \int_{t_{k-1} }^{t_k }
1 \cdot \int_{t_{k-1} }^{t_k } {\tilde {\omega }_i^2 (x,t)} \right)
}
 = \sum\limits_k {\int_{t_{k-1} }^{t_k } \int_\Omega {\tilde
{\omega }_i^2 (x,t)} } =\int_0^T \int_\Omega {\tilde
{\omega }_i^2 (x,t)} <+\infty \\
 \end{split}
\end{equation*}
\\

In the same way, we know from (\ref{eq2}) that the function
\[
\overline u_i (x,t)=\left\{ {\frac{1}{\Delta t_k }\int_{t_{k-1}
}^{t_k } { u_i (x,t)dt} ,\;\;t\in (t_{k-1} ,t_k )\subset (0,T)}
\right\}
\]
belongs to $L^2(0,T; H)$.
\\

Finally we will prove that $(\omega _1^\ast ,\omega _2^\ast ,\omega
_3^\ast )$ is a solution of the vorticity-velocity form of
Navier-Stokes equation (\ref{eq4}).
\\

Taking $\varphi _i \in C^\infty ((0,T)\times {\mathbb
R}^3),\;\;(i=1,2,3)$ with a period on $\Omega $, and
\[
\partial _{x_1 } \varphi _1 +\partial _{x_2 }
\varphi _2 +\partial _{x_3 } \varphi _3 =0
\]
we have
\begin{equation*}
\begin{split}
 &\sum\limits_{k=1}^N {\int_{t_{k-1} }^{t_k } {\int_\Omega {\varphi _1
(\partial _t \tilde {\omega }_1 \,+ \overline{\overline u}_1^k
\partial _{x_1 } \overline {\omega }_1^{k} + \overline{\overline u}_2^k \partial _{x_2 }
\overline {\omega }_1^{k} + \overline{\overline u}_3^k
\partial _{x_3 }
\overline {\omega }_1^{k} -} } } \\
\end{split}
\end{equation*}
\begin{equation*}
\begin{split}
 &\quad \quad \quad \quad \quad \quad \quad \quad \quad -\overline {\omega
}_1^{k} \partial _{x_1 } \overline{\overline u}_1^k -\overline
{\omega }_2^{k} \partial _{x_2 } \overline{\overline u}_1^k
-\overline {\omega }_3^{k} \partial _{x_3 } \overline{\overline
u}_1^k +\partial _{x_1 } q-\Delta
\tilde {\omega }_1 )=0 \\
 &\sum\limits_{k=1}^N {\int_{t_{k-1} }^{t_k } {\int_\Omega {\varphi _2
(\partial _t \tilde {\omega }_2 + \overline{\overline u}_1^k
\partial _{x_1 } \overline {\omega }_2^{k} + \overline{\overline u}_2^k \partial _{x_2 }
\overline {\omega }_2^{k} + \overline{\overline u}_3^k
\partial _{x_3 } \overline
{\omega }_2^{k} } } } - \\
 &\quad \quad \quad \quad \quad \quad \quad \quad \quad -\overline {\omega
}_1^{k} \partial _{x_1 } \overline{\overline u}_2^k -\overline
{\omega }_2^{k} \partial _{x_2 } \overline{\overline u}_2^k
-\overline {\omega }_3^{k} \partial _{x_3 } \overline{\overline
u}_2^k +\partial _{x_2 } q-\Delta
\tilde {\omega }_2 )=0 \\
 &\sum\limits_{k=1}^N {\int_{t_{k-1} }^{t_k } {\int_\Omega {\varphi _3
(\partial _t \tilde {\omega }_3 + \overline{\overline u}_1^k
\partial _{x_1 } \overline {\omega }_3^{k} + \overline{\overline u}_2^k \partial _{x_2 }
\overline {\omega }_3^{k} + \overline{\overline u}_3^k
\partial _{x_3 } \overline
{\omega }_3^{k} } } } - \\
 &\quad \quad \quad \quad \quad \quad \quad \quad \quad -\overline {\omega
}_1^{k} \partial _{x_1 } \overline{\overline u}_3^k -\overline
{\omega }_2^{k} \partial _{x_2 } \overline{\overline u}_3^k
-\overline {\omega }_3^{k} \partial _{x_3 } \overline{\overline
u}_3^k +\partial _{x_3 } q-\Delta
\tilde {\omega }_3 )=0 \\
 \end{split}
\end{equation*}
\\
Here $\tilde {\omega }_i \;(i=1,2,3)$ denote the collection of those
solutions of problem (5) defined on every $(t_{k-1},t_k)$.
Integrating by parts we get
\\
\begin{equation*}
\begin{split}
 &\sum\limits_{k=1}^N {\int_{t_{k-1} }^{t_k } {\int_\Omega {(\tilde {\omega
}_1 \partial _t \varphi _1 \,+\overline {\omega }_1^{k} ((
\overline{\overline u}_1^k
\partial _{x_1 } \varphi _1 +\varphi _1 \,\partial _{x_1 } \overline{\overline u}_1^k )+( \overline{\overline u}_2^k
\partial _{x_2 } \varphi _1 +\varphi _1 \,\partial _{x_2 } \overline{\overline u}_2^k )+( \overline{\overline u}_3^k \partial _{x_3 }
\varphi _1 +\varphi _1 \,\partial _{x_3 } \overline{\overline u}_3^k ))-} } } \\
 &\quad - \overline{\overline u}_1^k ((\overline {\omega }_1^{k} \partial _{x_1 } \varphi _1 +\varphi _1
\,\partial _{x_1 } \overline {\omega }_1^{k} )+(\overline {\omega
}_2^{k}
\partial _{x_2 } \varphi _1 +\varphi _1 \,\partial _{x_2 } \overline
{\omega }_2^{k} )+(\overline {\omega }_3^{k} \partial _{x_3 }
\varphi _1 +\varphi _1 \,\partial
_{x_3 } \overline {\omega }_3^{k} ))+ \\
 &\quad +q\partial _{x_1 } \varphi _1 +\tilde {\omega }_1 \Delta \varphi _1
)=\sum\limits_{k=1}^N {\int_\Omega {(\varphi _1 (x,t_k )\tilde {\omega }_1
(x,t_k )-\varphi _1 (x,t_{k-1} )\tilde {\omega }_1 (x,t_{k-1} ))} } \\
 &\sum\limits_{k=1}^N {\int_{t_{k-1} }^{t_k } {\int_\Omega {(\tilde {\omega
}_2 \partial _t \varphi _2 \,+\overline {\omega }_2^{k} ((
\overline{\overline u}_1^k
\partial _{x_1 } \varphi _2 +\varphi _2 \,\partial _{x_1 } \overline{\overline u}_1^k )+( \overline{\overline u}_2^k
\partial _{x_2 } \varphi _2 +\varphi _2 \,\partial _{x_2 } \overline{\overline u}_2^k )+( \overline{\overline u}_3^k \partial _{x_3 }
\varphi _2 +\varphi _2 \,\partial _{x_3 } \overline{\overline u}_3^k ))-} } } \\
 &\quad - \overline{\overline u}_2^k ((\overline {\omega }_1^{k} \partial _{x_1 } \varphi _2 +\varphi _2
\,\partial _{x_1 } \overline {\omega }_1^{k} )+(\overline {\omega
}_2^{k}
\partial _{x_2 } \varphi _2 +\varphi _2 \,\partial _{x_2 } \overline
{\omega }_2^{k} )+(\overline {\omega }_3^{k} \partial _{x_3 }
\varphi _2 +\varphi _2 \,\partial
_{x_3 } \overline {\omega }_3^{k} ))+ \\
 &\quad +q\partial _{x_2 } \varphi _2 +\tilde {\omega }_2 \Delta \varphi _2
)=\sum\limits_{k=1}^N {\int_\Omega {(\varphi _2 (x,t_k )\tilde {\omega }_2
(x,t_k )-\varphi _2 (x,t_{k-1} )\tilde {\omega }_2 (x,t_{k-1} ))} } \\
 &\sum\limits_{k=1}^N {\int_{t_{k-1} }^{t_k } {\int_\Omega {(\tilde {\omega
}_3 \partial _t \varphi _3 \,+\overline {\omega }_3^{k} ((
\overline{\overline u}_1^k
\partial _{x_1 } \varphi _3 +\varphi _3 \,\partial _{x_1 } \overline{\overline u}_1^k )+( \overline{\overline u}_2^k
\partial _{x_2 } \varphi _3 +\varphi _3 \,\partial _{x_2 } \overline{\overline u}_2^k )+( \overline{\overline u}_3^k \partial _{x_3 }
\varphi _3 +\varphi _3 \,\partial _{x_3 } \overline{\overline u}_3^k ))-} } } \\
 &\quad - \overline{\overline u}_3^k ((\overline {\omega }_1^{k} \partial _{x_1 } \varphi _3 +\varphi _3
\,\partial _{x_1 } \overline {\omega }_1^{k} )+(\overline {\omega
}_2^{k}
\partial _{x_2 } \varphi _3 +\varphi _3 \,\partial _{x_2 } \overline
{\omega }_2^{k} )+(\overline {\omega }_3^{k} \partial _{x_3 }
\varphi _3 +\varphi _3 \,\partial
_{x_3 } \overline {\omega }_3^{k} ))+ \\
 &\quad +q\partial _{x_3 } \varphi _3 +\tilde {\omega }_3 \Delta \varphi _3
)=\sum\limits_{k=1}^N {\int_\Omega {(\varphi _3 (x,t_k )\tilde {\omega }_3
(x,t_k )-\varphi _3 (x,t_{k-1} )\tilde {\omega }_3 (x,t_{k-1} ))} } \\
 \end{split}
\end{equation*}
\\

From section 2 we have the following conclusions:
\\

$\tilde {\omega }_i \to \omega _i^\ast $ in $L^2 (0,T;H)$ weakly,
and in $L^\infty (0,T;H)$ weak-star

$\overline {\omega }_i \to \omega _i^\ast $ in $L^2(0,T;H)$ weakly
\\
\\
as ${k^\prime} \to \infty $, or $\Delta t_k ^\prime \to 0$.
\\

In addition, for a certain solution $u$ of (\ref{eq1}), we can prove
due to (\ref{eq2}) and (\ref{eq3}) that
\\

$ \overline{\overline u}_i \to u_i $ in $L^2 (0,T;H )$ strongly
\\
\\
as $\varepsilon \to 0$ and ${k} \to \infty $, or $\Delta t_k \to 0$.
\\

In fact, we have
\begin{equation*}
\begin{split}
 &\left\| \,{\overline{\overline u} _i - \overline u_i } \right\|_{L^2(\Omega )} =
\left\| \;\; {\int\limits_{\left| y \right| \le \varepsilon } {J_\varepsilon
(y)\, [\,\overline u _i (x - y) - \overline u _i (x)\,]\, dy} } \right\|_{L^2(\Omega
)} \\
 &\le \int\limits_{\left| y \right| \le \varepsilon } {J_\varepsilon
(y)\left\| \, {\overline u _i (x - y) - \overline u _i (x)}
\right\|_{L^2(\Omega )} dy} \\
 &\le \mathop {\sup }\limits_{\left| y \right| \le
\varepsilon } \left\| \, {\overline u _i (x - y) - \overline u _i (x)}
\right\|_{L^2(\Omega )} \to 0 \\
 \end{split}
\end{equation*}
as $\varepsilon \to 0$. We can take $\varepsilon = O(\frac{1}{k}) $.
\\

Set $Q=(0,T)\times \overline {\Omega }$, $\Delta t=\mathop {\max
}\limits_k \{\Delta t_k \}$. $\forall \varepsilon >0$, and $u_i \in
L^2(0,T;L^2(\Omega ))$, there exists a $v_i \in C^\infty
(0,T;L^2(\Omega ))$ such that
\[
\left\| {u_i -v_i } \right\|_{L^2(Q)} <\varepsilon
\]
By means of the same partition as that for $\overline {u}_i $ to
construct $\overline {v}_i $, since there exists a constant $C>0$ such
that $\left\| {\partial _t v_i } \right\|_{L^2(\Omega )} \le C$, and
$\mathop {\max }\limits_t \,\,\left\| {\overline {v}_i -v_i }
\right\|_{L^2(\Omega )} \le C\;\Delta t$, it follows that
\[
\left\| {\overline {v}_i -v_i } \right\|_{L^2(Q)} =\left( {\int_0^T
{\left\| {\overline {v}_i -v_i } \right\|_{L^2(\Omega )}^2 } }
\right)^{1/2}\le CT^{1/2}\;\Delta t
\]
Thus
\[\overline {v}_i \to v_i \;\;\left( {L^\infty (0,T;L^2(\Omega ))}
\right),\;\;\;\; \mbox{as} \;\; \Delta t\to 0 \] Take $\Delta t$
such that $\left\| {\overline {v}_i -v_i } \right\|_{L^2(Q)} <\varepsilon
$. Moreover,
\begin{equation*}
\begin{split}
 &\int_0^T {\left\| {\overline {u}_i -\overline {v}_i } \right\|_{L^2(\Omega )}^2 }
=\sum\limits_{k=1}^N {\left\| {\frac{1}{\Delta t_k }\int_{t_{k-1}
}^{t_k }
{(u_i -v_i )} } \right\|} _{L^2(\Omega )}^2 \Delta t_k \\
 &\le \sum\limits_{k=1}^N {\left\| {\left( {\int_{t_{k-1} }^{t_k } {(u_i -v_i
)^2} } \right)^{1/2}} \right\|} _{L^2(\Omega )}^2 \le \int_0^T
{\left\| {u_i
-v_i } \right\|_{L^2(\Omega )}^2 } \\
 \end{split}
\end{equation*}
so that $\left\| {\overline {u}_i -\overline {v}_i } \right\|_{L^2(Q)} \le
\left\| {u_i -v_i } \right\|_{L^2(Q)} <\varepsilon $. Therefore,
\[
\left\| {\overline {u}_i -u_i } \right\|_{L^2(Q)} \le \left\| {u_i -v_i }
\right\|_{L^2(Q)} +\left\| {v_i -\overline {v}_i } \right\|_{L^2(Q)}
+\left\| {\overline {v}_i -\overline {u}_i } \right\|_{L^2(Q)} <3\varepsilon
\]
Hence as $\Delta t\to 0$, we have $\left\| {\overline {u}_i -u_i }
\right\|_{L^2(Q)} \to 0$.
\\

Finally we obtain
\[
\left\| \, {\overline{\overline u} _i - u_i } \right\|_{L^2(Q)} \le \left(
{\int_0^T {\left\| \, {\overline{\overline u} _i - \overline u _i }
\right\|_{L^2(\Omega )}^2 } } \right)^{1 / 2} +\; \left\| \, {\overline u _i -
u_i } \right\|_{L^2(Q)} \to 0
\]
as ${k} \to \infty $ or $\Delta t_k \to 0$.
\\

These convergence results enable us to pass the limit. That is,
\begin{equation*}
\begin{split}
 &\sum\limits_{k'} {\int_{t_{k'-1} }^{t_k' } {\int_\Omega {(\tilde {\omega
}_1 \partial _t \varphi _1 \,+\overline {\omega }_1^{k'} ( \overline{\overline u}_1^{k'}
\partial _{x_1 } \varphi _1 + \overline{\overline u}_2^{k'} \partial _{x_2 } \varphi _1 + \overline{\overline u}_3^{k'}
\partial _{x_3 } \varphi _1
)-} } } \\
 &\quad \quad \quad \quad \quad \quad - \overline{\overline u}_1^{k'} (\overline {\omega }_1^{k'} \partial
_{x_1 } \varphi _1 +\overline {\omega }_2^{k'} \partial _{x_2 } \varphi
_1 +\overline {\omega }_3^{k'} \partial _{x_3 } \varphi _1 )+q\partial
_{x_1 } \varphi _1
+\tilde {\omega }_1 \Delta \varphi _1 ) \\
 &\quad \quad \quad \quad \quad \quad =\int_\Omega {(\varphi _1 (x,T)\tilde
{\omega }_1 (x,T)-\varphi _1 (x,0)\tilde {\omega }_1 (x,0))} \\
 &\sum\limits_{k'} {\int_{t_{k'-1} }^{t_k' } {\int_\Omega {(\tilde {\omega
}_2 \partial _t \varphi _2 \,+\overline {\omega }_2^{k'} ( \overline{\overline u}_1^{k'}
\partial _{x_1 } \varphi _2 + \overline{\overline u}_2^{k'} \partial _{x_2 } \varphi _2 + \overline{\overline u}_3^{k'}
\partial _{x_3 } \varphi _2
)-} } } \\
 &\quad \quad \quad \quad \quad \quad - \overline{\overline u}_2^{k'} (\overline {\omega }_1^{k'} \partial
_{x_1 } \varphi _2 +\overline {\omega }_2^{k'} \partial _{x_2 } \varphi
_2 + \overline {\omega }_3^{k'} \partial _{x_3 } \varphi _2 )+q\partial
_{x_2 } \varphi _2
+ \tilde {\omega }_2 \Delta \varphi _2 ) \\
 &\quad \quad \quad \quad \quad \quad =\int_\Omega {(\varphi _2 (x,T)\tilde
{\omega }_2 (x,T)-\varphi _2 (x,0)\tilde {\omega }_2 (x,0))} \\
 &\sum\limits_{k'} {\int_{t_{k'-1} }^{t_k' } {\int_\Omega {(\tilde {\omega
}_3 \partial _t \varphi _3 \,+\overline {\omega }_3^{k'} ( \overline{\overline u}_1^{k'}
\partial _{x_1 } \varphi _3 + \overline{\overline u}_2^{k'} \partial _{x_2 } \varphi _3 + \overline{\overline u}_3^{k'}
\partial _{x_3 } \varphi _3
)-} } } \\
 &\quad \quad \quad \quad \quad \quad - \overline{\overline u}_3^{k'} (\overline {\omega }_1^{k'} \partial
_{x_1 } \varphi _3 +\overline {\omega }_2^{k'} \partial _{x_2 } \varphi
_3 + \overline {\omega }_3^{k'} \partial _{x_3 } \varphi _3 )+ q\partial
_{x_3 } \varphi _3
+\tilde {\omega }_3 \Delta \varphi _3 ) \\
 &\quad \quad \quad \quad \quad \quad =\int_\Omega {(\varphi _3 (x,T)\tilde
{\omega }_3 (x,T)-\varphi _3 (x,0)\tilde {\omega }_3 (x,0))} \\
 \end{split}
\end{equation*}
This is equivalent to
\begin{equation*}
\begin{split}
 &\int_0^T {\int_\Omega {\,\{(\omega _1^\ast \partial _t \varphi _1 \,+\omega
_2^\ast \partial _t \varphi _2 \,+\omega _3^\ast \partial _t \varphi _3 )+}
} \\
 &\quad \quad +(\omega _1^\ast \Delta \varphi _1 +\omega _2^\ast \Delta
\varphi _2 +\omega _3^\ast \Delta \varphi _3 )+ \\
 &\quad \quad +\omega _1^\ast (u_1 \partial _{x_1 } \varphi _1 +u_2 \partial
_{x_2 } \varphi _1 +u_3 \partial _{x_3 } \varphi _1 )+\omega _2^\ast (u_1
\partial _{x_1 } \varphi _2 +u_2 \partial _{x_2 } \varphi _2 +u_3 \partial
_{x_3 } \varphi _2 )+ \\
 &\quad \quad +\omega _3^\ast (u_1 \partial _{x_1 } \varphi _3 +u_2 \partial
_{x_2 } \varphi _3 +u_3 \partial _{x_3 } \varphi _3 ) \\
\end{split}
\end{equation*}
\begin{equation*}
\begin{split}
 &\quad \quad -u_1 (\omega _1^\ast \partial _{x_1 } \varphi _1 +\omega
_2^\ast \partial _{x_2 } \varphi _1 +\omega _3^\ast \partial _{x_3 } \varphi
_1 )-u_2 (\omega _1^\ast \partial _{x_1 } \varphi _2 +\omega _2^\ast
\partial _{x_2 } \varphi _2 +\omega _3^\ast \partial _{x_3 } \varphi _2 )-
\\
 &\quad \quad -u_3 (\omega _1^\ast \partial _{x_1 } \varphi _3 +\omega
_2^\ast \partial _{x_2 } \varphi _3 +\omega _3^\ast \partial _{x_3 } \varphi
_3 )\} \\
 &=\int_\Omega {\{(\varphi _1 (x,T)\omega _1^\ast (x,T)+\varphi _2
(x,T)\omega _2^\ast (x,T)+\varphi _3 (x,T)\omega _3^\ast (x,T))-} \\
 &\quad \quad \;\;-(\varphi _{10} (x)\omega _{10} (x)+\varphi _{20} (x)\omega
_{20} (x)+\varphi _{30} (x)\omega _{30} (x))\} \\
\end{split}
\end{equation*}
\\
Here we also have
\[
\omega _i^\ast (x,0)=\omega _{i0} (x),\quad \varphi _i (x,0)=\varphi _{i0}
(x),\quad i=1,2,3
\]
Hence we know that there exists some $\omega _i^\ast $ which belongs
to $L^\infty (0,T;L^2(\Omega ))$ and is a Leray-Hopf weak solution
of (\ref{eq4}).
\\

Note that a weak formulation of the following equations:
\begin{equation*}
\begin{split}
 &\;\; \omega = \mbox{curl}\,u   \\
 &\int_0^{T} {\int_\Omega {\;\varphi \cdot [\partial _t \omega + (u \cdot
\nabla )\omega - (\omega \cdot \nabla )u - \nu\,\Delta \omega } } ]
= 0
\\
\end{split}
\end{equation*}
are equivalent to
\[
\int_0^{T} {\int_\Omega {\;\tilde {\varphi } \cdot [\partial _t u +
(u \cdot \nabla )u + \nabla p - \nu\,\Delta u} } ] = 0
\]
for any $\varphi \in C^\infty ((0,{T})\times {\rm \mathbb{R}}^3)$
with a period on $\Omega $, and $\tilde {\varphi } =
\mbox{curl}\varphi $, in some distribution sense.
\\
\\
\\

\textbf{5. Nonexistence}

Consider a series
\[
\zeta (s) = \sum\limits_{n = 1}^\infty {\frac{1}{n^s}}
\]
where $s$ is a complex variable, $s = \sigma + \mbox{i}t$. Its analytic continuation is called Riemann zeta-function. The existence and infinity of its nontrivial zeros on the critical line $\sigma = \frac{1}{2}$ are already known.
\\

Let
\[
\xi (s)\;\; = \frac{1}{2}s(s - 1)\;\pi ^{ - \frac{1}{2}s}\;\Gamma \left(
{\frac{1}{2}s} \right)\;\zeta (s)
\]
which is an entire function, and all zeros of $\xi (s)$ coincide with the nontrivial zeros of $\zeta (s)$.
\\

Writing
\[
\Xi (z) = \xi \left( {\frac{1}{2} + \mbox{i}z} \right),\quad z = t -
\mbox{i}\tilde {\sigma },\quad \tilde {\sigma } = \sigma - \frac{1}{2}
\]
then, as well known, we have an analytic expression of $\Xi (s)$ without those poles
and trivial zeros of $\zeta (s)$ as follows :
\begin{equation*}
\begin{split}
\qquad \qquad \qquad \qquad    \Xi (z) = 2\int_0^\infty {\Phi (t)\cos (zt)dt}    \qquad \qquad \qquad \qquad \qquad \qquad  (*)
\end{split}
\end{equation*}
where
\[
\Phi (t) = 2\;\sum\limits_{n = 1}^\infty {(\,2n^4\pi ^2e^{\frac{9}{2}t} -
3n^2\pi e^{\frac{5}{2}t}\,)\;e^{ - n^2\pi \,e^{2t}}}
\]
\\

Let $z = y - \mbox{i}x$, then
\begin{equation*}
\begin{split}
 &\cos z = \frac{1}{2}(e^{\mbox{i}\,z} + e^{ - \mbox{i}\,z}) =
\frac{1}{2}(e^{x + \mbox{i}\,y} + e^{ - x - \mbox{i}\,y}) \\
 &\quad \quad = \frac{1}{2}\;[\,e^x(\cos y + \mbox{i}\sin y) + e^{ - x}(\cos
y - \mbox{i}\sin y)\,] \\
 &\quad \quad = \frac{e^x + e^{ - x}}{2}\cos y + \;\mbox{i}\;\frac{e^x - e^{
- x}}{2}\sin y \\
 &\quad \quad = \cosh x\cos y + \mbox{i}\sinh x\sin y \\
\end{split}
\end{equation*}
where $\cosh x$ and $\sinh x$ are two hyperbolic functions.
\\

From (*) we know that
\begin{equation*}
\begin{split}
 &\Xi (z) = \xi \left( {\frac{1}{2} + \mbox{i}z} \right) = \xi \left(
{\frac{1}{2} + x + \mbox{i}y} \right) \\
 &\quad \quad = 2\;\int_0^\infty {\Phi (t)\cos (zt)\,dt} \\
\end{split}
\end{equation*}
where
\[
\Phi (t) = 2\;\sum\limits_{n = 1}^\infty {(\,2n^4\pi ^2e^{\frac{9}{2}t} -
3n^2\pi \,e^{\frac{5}{2}t}\,)\;e^{ - n^2\pi \,e^{2t}}}
\]
It follows that
\begin{equation*}
\begin{split}
 &\Xi (z) = 2\;\int_0^\infty {\Phi (t)\;\{\,\cosh (xt)\cos (yt) +
\mbox{i}\sinh (xt)\sin (yt)\,\}} \,dt \\
 &\quad \quad = 2\;\int_0^\infty {\Phi (t)\;\cosh (xt)\cos (yt)\,dt} +
\mbox{i}\;2\;\int_0^\infty {\Phi (t)\,\sinh (xt)\sin (yt)\,dt} \\
\end{split}
\end{equation*}
Thus $\Xi (z)$ is divided into the real and imaginary parts as
\begin{equation*}
\begin{split}
 &f(x,y) = 2\,\int_0^\infty {\Phi (t)\,\cosh (xt)\cos (yt)\,dt} \\
 &g(x,y) = 2\,\int_0^\infty {\Phi (t)\,\sinh (xt)\sin (yt)\,dt} \\
\end{split}
\end{equation*}
then $\Xi (z) = 0$ is equivalent to the following system of equations :
\[
f(x,y) = 0,\quad g(x,y) = 0
\]

For any fixed $x_0 > 0$, we will consider a property of zeros of $f(x_0 ,y)$ and $g(x_0,y)$ in this section.
\\

\textbf{Lemma 1. }  For a fixed $x_0 > 0$, $f(x_0 ,y)$ or $g(x_0 ,y)$ with respect to
$y$ is not of any zero of second or higher even order.
\\

In fact, let the arbitrary $\delta > 0$ be small enough, if $f(x_0 ,y_0 ) =
0$, then
\begin{equation*}
\begin{split}
 &\qquad f(x_0 ,y_0 - \delta ) \cdot f(x_0 ,y_0 + \delta ) \\
 &= \int_0^\infty {\Phi (t)\cosh (x_0 t)\cos (y_0 - \delta )t\,dt} \;\, \cdot
\;\,\int_0^\infty {\Phi (t)\cosh (x_0 t)\cos (y_0 + \delta )t\,dt} \\
 &= \int_0^\infty {\Phi (t)\cosh (x_0 t)\;[\,\cos (y_0 t)\cos (\delta t) +
\,\sin (y_0 t)\sin (\delta t)\,]\,dt} \;\,\times \\
 &\quad \;\times \int_0^\infty {\Phi (t)\cosh (x_0 t)\;[\,\cos (y_0 t)\cos
(\delta t) - \,\sin (y_0 t)\sin (\delta t)\,]\,dt} \\
 &= \left( {\int_0^\infty {\Phi (t)\cosh (x_0 t)\cos (y_0 t)\cos (\delta
t)\,dt} } \right)^2 - \left( {\int_0^\infty {\Phi (t)\cosh (x_0 t)\sin (y_0
t)\sin (\delta t)\,dt} } \right)^2 \\
\end{split}
\end{equation*}

An application of the following asymptotic expansion
\begin{equation*}
\begin{split}
 &\qquad \qquad  \sin x = x - \frac{1}{3!}\theta _1 (x)x^3,\quad \quad \quad \quad \quad
\left| {\theta _1 (x)} \right| \le 1 \\
 &\qquad \qquad  1 - \cos x = \frac{1}{2!}x^2 - \frac{1}{4!}\theta _2 (x)x^4,\quad \quad
\left| {\theta _2 (x)} \right| \le 1          \qquad \qquad \qquad \qquad      (**)    \\
\end{split}
\end{equation*}
leads to
\begin{equation*}
\begin{split}
 &\qquad f(x_0 ,y_0 - \delta )\; \cdot \;f(x_0 ,y_0 + \delta ) \\
 &= \left( {\int_0^\infty {\Phi (t)\cosh (x_0 t)\cos (y_0 t)(1 - \cos \delta
t)\,dt} } \right)^2 - \left( {\int_0^\infty {\Phi (t)\cosh (x_0 t)\sin (y_0
t)\sin (\delta t)\,dt} } \right)^2 \\
 &= \left( {\int_0^\infty {\Phi (t)\cosh (x_0 t)\cos (y_0
t)\;[\,\frac{1}{2!}(\delta t)^2 - \frac{1}{4!}\,\theta _2 (\delta
t)\,(\delta t)^4\,]\,dt} } \right)^2 - \\
 &\quad \quad - \left( {\int_0^\infty {\Phi (t)\cosh (x_0 t)\sin (y_0
t)\;[\,(\delta t) - \frac{1}{3!}\,\theta _1 (\delta t)\,(\delta t)^3\,]\,dt}
} \right)^2 \\
 &= \delta ^4\left( {\int_0^\infty {t^2\,\Phi (t)\cosh (x_0 t)\cos (y_0
t)\;[\,\frac{1}{2} - \frac{1}{24}\,\theta _2 (\delta t)\,(\delta
t)^2\,]\,dt} } \right)^2 - \\
 &\quad \quad - \delta ^2\left( {\int_0^\infty {t\,\Phi (t)\cosh (x_0 t)\sin
(y_0 t)\;[\,1 - \frac{1}{6}\,\theta _1 (\delta t)\,(\delta t)^2\,]\,dt} }
\right)^2 \\
 &< 0 \\
\end{split}
\end{equation*}
Thus, as $f(x_0 ,y_0 - \delta ) > 0,\;\;f(x_0 ,y_0 + \delta ) < 0$, or vice
versa. The conclusion is obtained.
\\

Similarly, using (**) we also have
\begin{equation*}
\begin{split}
 &\qquad g(x_0 ,y_0 - \delta )\; \cdot \;g(x_0 ,y_0 + \delta ) \\
 &= \int_0^\infty {\Phi (t)\sinh (x_0 t)\sin (y_0 - \delta )t\,dt} \;\, \cdot
\,\;\int_0^\infty {\Phi (t)\sinh (x_0 t)\sin (y_0 + \delta )t\,dt} \\
\end{split}
\end{equation*}
\begin{equation*}
\begin{split}
 &= \int_0^\infty {\Phi (t)\sinh (x_0 t)\;[\,\sin (y_0 t)\cos (\delta t) -
\cos (y_0 t)\sin (\delta t)\,]\,\,dt} \;\;\times \\
 &\quad \quad \times \int_0^\infty {\Phi (t)\sinh (x_0 t)\;[\,\sin (yt)\cos
(\delta t) + \cos (y_0 t)\sin (\delta t)\,]\,\,dt} \\
 &= \left( {\int_0^\infty {\Phi (t)\sinh (x_0 t)\;\sin (y_0 t)\cos (\delta
t)\,dt} } \right)^2 - \left( {\int_0^\infty {\Phi (t)\sinh (x_0 t)\;\cos
(y_0 t)\sin (\delta t)\,dt} } \right)^2 \\
 &= \left( {\int_0^\infty {\Phi (t)\sinh (x_0 t)\;\sin (y_0 t)\,(1 - \cos
\delta t)\,dt} } \right)^2 - \left( {\int_0^\infty {\Phi (t)\sinh (x_0
t)\;\cos (y_0 t)\sin (\delta t)\,dt} } \right)^2 \\
 &= \left( {\int_0^\infty {\Phi (t)\sinh (x_0 t)\;\sin (y_0
t)\;[\,\frac{1}{2!}(\delta t)^2 - \frac{1}{4!}\theta _2 (\delta t)\,(\delta
t)^4]\,dt} } \right)^2 - \\
 &\quad \quad - \left( {\int_0^\infty {\Phi (t)\sinh (x_0 t)\;\cos (y_0
t)\;[\,(\delta t) - \frac{1}{3!}\theta _1 (\delta t)\,(\delta t)^3]\,dt} }
\right)^2 \\
 &= \delta ^4\;\left( {\int_0^\infty {t^2\Phi (t)\sinh (x_0 t)\;\sin (y_0
t)\;[\,\frac{1}{2} - \frac{1}{24}\theta _2 (\delta t)\,(\delta t)^2]\,dt} }
\right)^2 - \\
 &\quad \quad - \delta ^2\;\left( {\int_0^\infty {t\,\Phi (t)\sinh (x_0
t)\;\cos (y_0 t)\;[\,1 - \frac{1}{6}\theta _1 (\delta t)\,(\delta t)^2]\,dt}
} \right)^2 \\
 &< 0 \\
\end{split}
\end{equation*}
Thus, the same conclusion is obtained.
\\

\textbf{Lemma 2. }  For a fixed $x_0 > 0$, $f(x_0 ,y)$ or $g(x_0 ,y)$ with respect to
$y$ is not of any zero of third or higher odd order.
\\

We assume for sake of contradiction that this conclusion fails, then ${f}'_y
(x_0 ,y_0 ) = 0$, and consider further that
\begin{equation*}
\begin{split}
 &\qquad {f}'_y (x_0 ,y_0 - \delta )\; \cdot \;{f}'_y (x_0 ,y_0 + \delta ) \\
 &= \left\{ { - \int_0^\infty {t\,\Phi (t)\cosh (x_0 t)\sin (y_0 - \delta
)t\,dt} } \right\}\; \cdot \;\left\{ { - \int_0^\infty {t\,\Phi (t)\cosh
(x_0 t)\sin (y_0 + \delta )t\,dt} } \right\} \\
 &= \int_0^\infty {t\,\Phi (t)\cosh (x_0 t)\;[\,\sin (y_0 t)\cos (\delta t) -
\cos (y_0 t)\sin (\delta t)\,]\,dt} \;\;\times \\
 &\quad \quad \times \int_0^\infty {t\,\Phi (t)\cosh (x_0 t)\;[\,\sin (y_0
t)\cos (\delta t) + \cos (y_0 t)\sin (\delta t)\,]\,dt} \\
 &= \left( {\int_0^\infty {t\,\Phi (t)\cosh (x_0 t)\;\sin (y_0 t)\cos (\delta
t)\,dt} } \right)^2 - \left( {\int_0^\infty {t\,\Phi (t)\cosh (x_0 t)\;\cos
(y_0 t)\sin (\delta t)\,dt} } \right)^2 \\
 &= \left( {\int_0^\infty {t\,\Phi (t)\cosh (x_0 t)\;\sin (y_0 t)\,(1 - \cos
\delta t)\,dt} } \right)^2 - \left( {\int_0^\infty {t\,\Phi (t)\cosh (x_0
t)\;\cos (y_0 t)\sin (\delta t)\,dt} } \right)^2 \\
 &= \left( {\int_0^\infty {t\,\Phi (t)\cosh (x_0 t)\;\sin (y_0
t)\;[\,\frac{1}{2!}(\delta t)^2 - \frac{1}{4!}\theta _2 (\delta t)\,(\delta
t)^4]\,dt} } \right)^2 - \\
\end{split}
\end{equation*}
\begin{equation*}
\begin{split}
 &\quad \quad - \left( {\int_0^\infty {t\,\Phi (t)\cosh (x_0 t)\;\cos (y_0
t)\;[\,(\delta t) - \frac{1}{3!}\theta _1 (\delta t)\,(\delta t)^3]\,dt} }
\right)^2 \\
 &= \delta ^4\left( {\int_0^\infty {t^3\,\Phi (t)\cosh (x_0 t)\;\sin (y_0
t)\;[\,\frac{1}{2} - \frac{1}{24}\theta _2 (\delta t)\,(\delta t)^2]\,dt} }
\right)^2 - \\
 &\quad \quad - \delta ^2\left( {\int_0^\infty {t^2\,\Phi (t)\cosh (x_0
t)\;\cos (y_0 t)\;[\,1 - \frac{1}{6}\theta _1 (\delta t)\,(\delta t)^2]\,dt}
} \right)^2 \\
 &< 0 \\
\end{split}
\end{equation*}
It follows that ${f}'_y (x_0 ,y_0 - \delta ) > 0$, if ${f}'_y (x_0 ,y_0 +
\delta ) < 0$, or vice versa. This is inconsistent with the previous
assumption that ${f}'_y (x_0 ,y)$ is of a zero of second or higher even
order.
\\

Similarly, we take proof by contradiction for $g(x_0 ,y)$ to assume that the
conclusion fails, then ${g}'_y (x_0 ,y_0 ) = 0$
\begin{equation*}
\begin{split}
 &\qquad {g}'_y (x_0 ,y_0 - \delta )\; \cdot \;{g}'_y (x_0 ,y_0 + \delta ) \\
 &= \int_0^\infty {t\,\Phi (t)\sinh (x_0 t)\cos (y_0 - \delta )t\,dt} \;\;
\cdot \;\;\int_0^\infty {t\,\Phi (t)\sinh (x_0 t)\cos (y_0 + \delta )t\,dt}
\\
 &= \int_0^\infty {t\,\Phi (t)\sinh (x_0 t)\;[\,\cos (y_0 t)\cos (\delta t) +
\sin (y_0 t)\sin (\delta t)\,]\,dt} \;\;\times \\
 &\quad \quad \times \int_0^\infty {t\,\Phi (t)\sinh (x_0 t)\;[\,\cos (y_0
t)\cos (\delta t) - \sin (y_0 t)\sin (\delta t)\,]\,dt} \\
 &= \left( {\int_0^\infty {t\,\Phi (t)\sinh (x_0 t)\;\cos (y_0 t)\cos (\delta
t)\,dt} } \right)^2 - \left( {\int_0^\infty {t\,\Phi (t)\sinh (x_0 t)\;\sin
(y_0 t)\sin (\delta t)\,dt} } \right)^2 \\
 &= \left( {\int_0^\infty {t\,\Phi (t)\sinh (x_0 t)\;\cos (y_0 t)\;(1 - \cos
\delta t)\,dt} } \right)^2 - \left( {\int_0^\infty {t\,\Phi (t)\sinh (x_0
t)\;\sin (y_0 t)\sin (\delta t)\,dt} } \right)^2 \\
 &= \left( {\int_0^\infty {t\,\Phi (t)\sinh (x_0 t)\;\cos (y_0
t)\;[\,\frac{1}{2!}(\delta t)^2 - \frac{1}{4!}\theta _2 (\delta t)(\delta
t)^4]\,dt} } \right)^2 - \\
 &\quad \quad - \left( {\int_0^\infty {t\,\Phi (t)\sinh (x_0 t)\;\sin (y_0
t)\;[\,(\delta t) - \frac{1}{3!}\theta _1 (\delta t)(\delta t)^3]\,dt} }
\right)^2 \\
 &= \delta ^4\;\left( {\int_0^\infty {t^3\,\Phi (t)\sinh (x_0 t)\;\cos (y_0
t)\;[\,\frac{1}{2} - \frac{1}{24}\theta _2 (\delta t)(\delta t)^2]\,dt} }
\right)^2 - \\
 &\quad \quad - \delta ^2\left( {\int_0^\infty {t^2\,\Phi (t)\sinh (x_0
t)\;\sin (y_0 t)\;[\,1 - \frac{1}{6}\theta _1 (\delta t)(\delta t)^2]\,dt} }
\right)^2 \\
 &< 0 \\
\end{split}
\end{equation*}
It follows that ${g}'_y (x_0 ,y)$ changes sign near the point $y_0 $. This
is also inconsistent with the assumption that ${g}'_y (x_0 ,y)$ is of a zero
of second or higher even order.
\\

It is impossible for an analytic function to be of an zero of fractional order, otherwise a singularity of some its derivative will appear. Combining Lemma 1 and 2, we have

\textbf{Theorem 1. } For any fixed $x_0 > 0$, all zeros of $f(x_0 ,y)$ and $g(x_0
,y)$ with respect to $y$ are simple.
\\

In the following, we deal with the nonexistence of any zero of $\Xi (z)$
standing outside the critical line $x = 0$ or $\sigma = \frac{1}{2}$. Since
$\Xi (z)$ is an even function, we only need considering the half-plane $x > 0$.
\\

We can prove that the system of equations $f(x,y) = 0,\;\,g(x,y) = 0$ has no real
solution as $x > 0$. That is, when $x > 0$, $f(x,y)$ and $g(x,y)$ can not
vanish at the same time.
\\

In fact, given any a point $(x_0 ,y_0 )$ on the half-plane $x > 0$, if $g(x_0 ,y_0 )
\ne 0$, then the system of equations has no solution. Otherwise, $g(x_0 ,y_0
) = 0$, that is, $(x_0 ,y_0 )$ satisfies some equation $g(x,y) = 0$. Since
$g(x,y)$ is continuous and differentiable, then there exists a neighborhood
\[
N_{\delta _1}(x_0, y_0) \;= \{ (x,y) : \sqrt {(x - x_0 )^2 + (y - y_0 )^2} < \delta _1 \}
\]
such that ${g}'_y (x,y) \ne 0$, according to Theorem 1. From
existence theorem of implicit function it follows that the equation $g(x,y)
= 0$ decides a continuous and smooth curve $y = \phi (x)$ in the
neighborhood. If $f(x_0 ,y_0 ) = 0$ at the same time, then ${f}'_y (x_0 ,y_0
) \ne 0$ according to Theorem 1. By means of the Cauchy-Riemann
equation, we have
\[
{g}'_x (x_0 ,y_0 ) = - {f}'_y (x_0 ,y_0 ) \ne 0
\]
Since $ {g}'_x (x,y) $ is also continuous and differentiable, then there exists another neighborhood with the same center
\[
N_{\delta _2}(x_0, y_0) \;= \{ (x,y) : \sqrt {(x - x_0 )^2 + (y - y_0 )^2} < \delta _2 \le \delta _1   \}
\]
such that ${g}'_x (x,y) \ne 0$. Let the point $ (x_1, y_1) \in N_{\delta _2}(x_0, y_0) $ be on the curve $y = \phi (x)$ and $ x_1 \ne x_0 $, as $ \delta _2 $ is small enough, we have
\[
 g(x_1,y_1 ) = {g}'_x (x_0 ,y_1 )\,(x_1 - x_0 )\; + \; O((x_1 - x_0 )^2)   \;\ne 0   \\
\]
This leads to a contradiction.
\\

As a conclusion, we obtain

\textbf{Theorem 2. } For any $z = y - \mbox{i}x,\;\,x > 0$, the system of
equations $f(x,y) = 0$ and $g(x,y) = 0$ has no real solution. That is,
\[
\Xi (z) = \xi \left( {\frac{1}{2} + x + \mbox{i}y} \right) \ne 0
\]
\\
\\
\\

\textbf{6. Regularity}

We can still use Galerkin procedure as in section 3. Since $V$ is
separable there exists a sequence of linearly independent elements
$w_{i1} ,\cdots ,w_{im} ,\cdots $ which is total in $V$. For each
$m$ we define an approximate solution $u_{im} $ of (1) as follows:
\[
u_{im} =\sum\limits_{j=1}^m {g_{ij} (t)w_{ij} }
\]
and
\begin{equation}
\label{eq14}
\begin{split}
 &\int_\Omega {\partial _t u_{1m} w_{1j} } +\int_\Omega {(u_{1m} \partial
_{x_1 } u_{1m} +u_{2m} \partial _{x_2 } u_{1m} +u_{3m} \partial
_{x_3 } u_{1m} )} \,w_{1j} +\int_\Omega {\partial _{x_1 } p\,w_{1j}
} =\int_\Omega
{\Delta u_{1m} \,w_{1j} } \\
 &\int_\Omega {\partial _t u_{2m} w_{2j} } +\int_\Omega {(u_{1m} \partial
_{x_1 } u_{2m} +u_{2m} \partial _{x_2 } u_{2m} +u_{3m} \partial
_{x_3 } u_{2m} )} \,w_{2j} +\int_\Omega {\partial _{x_2 } p\,w_{2j}
} =\int_\Omega
{\Delta u_{2m} \,w_{2j} } \\
 &\int_\Omega {\partial _t u_{3m} w_{3j} } +\int_\Omega {(u_{1m} \partial
_{x_1 } u_{3m} +u_{2m} \partial _{x_2 } u_{3m} +u_{3m} \partial
_{x_3 } u_{3m} )} \,w_{3j} +\int_\Omega {\partial _{x_3 } p\,w_{3j}
} =\int_\Omega
{\Delta u_{3m} \,w_{3j} } \\
 &\quad \quad u_{im} (0)=u_{i0}^m ,\quad \quad j=1,\cdots ,m \\
 \end{split}
\end{equation}
where $u_{i0}^m $ is the orthogonal projection in $H$ of $u_{i0} $
on the space spanned by $w_{i1} ,\cdots ,w_{im} $.
\\

We now are allowed to differentiate (\ref{eq14}) in the $t$, we get
\begin{equation*}
\begin{split}
 &\int_\Omega {\partial _t^2 u_{1m} w_{1j} } +\int_\Omega {(\partial _t
u_{1m} \partial _{x_1 } u_{1m} +\partial _t u_{2m} \partial _{x_2 }
u_{1m}
+\partial _t u_{3m} \partial _{x_3 } u_{1m} )} \,w_{1j} + \\
 &\quad +\int_\Omega {(u_{1m} \partial _{x_1 } \partial _t u_{1m} +u_{2m}
\partial _{x_2 } \partial _t u_{1m} +u_{3m} \partial _{x_3 } \partial _t
u_{1m} )} \,w_{1j} +\int_\Omega {\partial _{x_1 } \partial _t
p\,w_{1j} }
=\int_\Omega {\Delta \partial _t u_{1m} \,w_{1j} } \\
\end{split}
\end{equation*}
\begin{equation}
\label{eq15}
\begin{split}
 &\int_\Omega {\partial _t^2 u_{2m} w_{2j} } +\int_\Omega {(\partial _t
u_{1m} \partial _{x_1 } u_{2m} +\partial _t u_{2m} \partial _{x_2 }
u_{2m}
+\partial _t u_{3m} \partial _{x_3 } u_{2m} )} \,w_{2j} + \\
 &\quad +\int_\Omega {(u_{1m} \partial _{x_1 } \partial _t u_{2m} +u_{2m}
\partial _{x_2 } \partial _t u_{2m} +u_{3m} \partial _{x_3 } \partial _t
u_{2m} )} \,w_{2j} +\int_\Omega {\partial _{x_2 } \partial _t
p\,w_{2j} }
=\int_\Omega {\Delta \partial _t u_{2m} \,w_{2j} } \\
 &\int_\Omega {\partial _t^2 u_{3m} w_{3j} } +\int_\Omega {(\partial _t
u_{1m} \partial _{x_1 } u_{3m} +\partial _t u_{2m} \partial _{x_2 }
u_{3m}
+\partial _t u_{3m} \partial _{x_3 } u_{3m} )} \,w_{3j} + \\
 &\quad +\int_\Omega {(u_{1m} \partial _{x_1 } \partial _t u_{3m} +u_{2m}
\partial _{x_2 } \partial _t u_{3m} +u_{3m} \partial _{x_3 } \partial _t
u_{3m} )} \,w_{3j} +\int_\Omega {\partial _{x_3 } \partial _t
p\,w_{3j} }
=\int_\Omega {\Delta \partial _t u_{3m} \,w_{3j} } \\
 &\qquad \qquad \qquad \qquad \qquad \qquad \qquad \qquad  j=1,\cdots ,m \\
 \end{split}
\end{equation}
We multiply (15) by ${g}'_{ij} (t)$ and add the resulting equations
for $j=1,\cdots ,m$, we find
\begin{equation*}
\begin{split}
 &\frac{1}{2}\partial _t \int_\Omega {(\partial _t u_{1m} )^2} +\int_\Omega
{\partial _t u_{1m} (\partial _t u_{1m} \partial _{x_1 } u_{1m}
+\partial _t u_{2m}
\partial _{x_2 } u_{1m} +\partial _t u_{3m} \partial _{x_3 } u_{1m} )} + \\
 &\quad \quad +\int_\Omega {\partial _t u_{1m} (u_{1m} \partial _{x_1 } \partial _t
u_{1m} +u_{2m} \partial _{x_2 } \partial _t u_{1m} +u_{3m} \partial
_{x_3 }
\partial _t u_{1m} )} +\int_\Omega {\partial _t u_{1m} \partial _{x_1 }
\partial _t p}
=\int_\Omega {\partial _t u_{1m} \,\Delta \partial _t u_{1m} } \\
 &\frac{1}{2}\partial _t \int_\Omega {(\partial _t u_{2m} )^2} +\int_\Omega
{\partial _t u_{2m} (\partial _t u_{1m} \partial _{x_1 } u_{2m}
+\partial _t u_{2m}
\partial _{x_2 } u_{2m} +\partial _t u_{3m} \partial _{x_3 } u_{2m} )} + \\
 &\quad \quad +\int_\Omega {\partial _t u_{2m} (u_{1m} \partial _{x_1 } \partial _t
u_{2m} +u_{2m} \partial _{x_2 } \partial _t u_{2m} +u_{3m} \partial
_{x_3 }
\partial _t u_{2m} )} +\int_\Omega {\partial _t u_{2m} \partial _{x_2 }
\partial _t p}
=\int_\Omega {\partial _t u_{2m} \,\Delta \partial _t u_{2m} } \\
 &\frac{1}{2}\partial _t \int_\Omega {(\partial _t u_{3m} )^2} +\int_\Omega
{\partial _t u_{3m} (\partial _t u_{1m} \partial _{x_1 } u_{3m}
+\partial _t u_{2m}
\partial _{x_2 } u_{3m} +\partial _t u_{3m} \partial _{x_3 } u_{3m} )} + \\
 &\quad \quad +\int_\Omega {\partial _t u_{3m} (u_{1m} \partial _{x_1 } \partial _t
u_{3m} +u_{2m} \partial _{x_2 } \partial _t u_{3m} +u_{3m} \partial
_{x_3 }
\partial _t u_{3m} )} +\int_\Omega {\partial _t u_{3m} \partial _{x_3 }
\partial _t p}
=\int_\Omega {\partial _t u_{3m} \,\Delta \partial _t u_{3m} } \\
 \end{split}
\end{equation*}
and
\[
\int_\Omega {(\partial _t u_{1m} \partial _{x_1 } \partial _t
p+\partial _t u_{2m}
\partial _{x_2 } \partial _t p+\partial _t u_{3m} \partial _{x_3 } \partial _t
p)} =-\int_\Omega {\partial _t p\,\partial _t (\partial _{x_1 }
u_{1m} +\partial _{x_2 } u_{2m} +\partial _{x_3 } u_{3m} )} =0
\]

Since
\begin{equation*}
\begin{split}
 &\int_\Omega {\partial _t u_{im} (u_{1m} \partial _{x_1 } \partial _t u_{im}
 +u_{2m}
\partial _{x_2 } \partial _t u_{im} +u_{3m} \partial _{x_3 } \partial _t u_{im} )} =
\\
 &\quad =\frac{1}{2}\int_\Omega {(u_{1m} \partial _{x_1 } (\partial _t
 u_{im}
)^2+u_{2m} \partial _{x_2 } (\partial _t u_{im} )^2+u_{3m} \partial
_{x_3 } (\partial
_t u_{im} )^2)} \\
 &\quad =-\frac{1}{2}\int_\Omega {(\partial _t u_{im} )^2(\partial _{x_1 }
 u_{1m}
+\partial _{x_2 } u_{2m} +\partial _{x_3 } u_{3m} )} =0 \\
 \end{split}
\end{equation*}
and
\begin{equation*}
\begin{split}
 &\int_\Omega {\partial _t u_{im} \,\Delta \partial _t u_{im} } =\int_\Omega
{\partial _t u_{im} (\partial _{x_1 }^2 \partial _t u_{im} +\partial
_{x_2 }^2
\partial _t u_{im} +\partial _{x_3 }^2 \partial _t u_{im} )} = \\
 &=-\int_\Omega {((\partial _{x_1 } \partial _t u_{im} )^2+(\partial _{x_2 }
\partial _t u_{im} )^2+(\partial _{x_3 } \partial _t u_{im} )^2)} ,\quad \quad
i=1,2,3 \\
 \end{split}
\end{equation*}
then
\begin{equation*}
\begin{split}
 &\frac{1}{2}\partial _t \int_\Omega {((\partial _t u_{1m} )^2+(\partial _t
 u_{2m}
)^2+(\partial _t u_{3m} )^2)} + \\
 &\qquad\quad +\left\| {\nabla \partial _t u_{1m} } \right\|_{L^2(\Omega )}^2 +\left\|
{\nabla \partial _t u_{2m} } \right\|_{L^2(\Omega )}^2 +\left\|
{\nabla
\partial _t u_{3m} } \right\|_{L^2(\Omega )}^2 \\
 &\le \left\| {\partial _t u_{1m} } \right\|_{L^4(\Omega )} \left( {\left\|
{\partial _t u_{1m} } \right\|_{L^4(\Omega )} \left\| {\partial
_{x_1 } u_{1m} } \right\|_{L^2(\Omega )} +\left\| {\partial _t
u_{2m} } \right\|_{L^4(\Omega )} \left\| {\partial _{x_2 } u_{1m} }
\right\|_{L^2(\Omega )} \;+ } \right. \\
 & \qquad\qquad\qquad\qquad\qquad\qquad  + \left. {\left\| {\partial _t u_{3m} }
\right\|_{L^4(\Omega )} \left\| {\partial _{x_3 } u_{1m} }
\right\|_{L^2(\Omega )} } \right) \\
 &+\left\| {\partial _t u_{2m} } \right\|_{L^4(\Omega )} \left( {\left\|
{\partial _t u_{1m} } \right\|_{L^4(\Omega )} \left\| {\partial
_{x_1 } u_{2m} } \right\|_{L^2(\Omega )} +\left\| {\partial _t
u_{2m} } \right\|_{L^4(\Omega )} \left\| {\partial _{x_2 } u_{2m} }
\right\|_{L^2(\Omega )} \;+ } \right. \\
 &\qquad\qquad\qquad\qquad\qquad\qquad + \left. {\left\| {\partial _t u_{3m} }
\right\|_{L^4(\Omega )} \left\| {\partial _{x_3 } u_{2m} }
\right\|_{L^2(\Omega )} } \right) \\
\end{split}
\end{equation*}
\begin{equation*}
\begin{split}
 &+\left\| {\partial _t u_{3m} } \right\|_{L^4(\Omega )} \left( {\left\|
{\partial _t u_{1m} } \right\|_{L^4(\Omega )} \left\| {\partial
_{x_1 } u_{3m} } \right\|_{L^2(\Omega )} +\left\| {\partial _t
u_{2m} } \right\|_{L^4(\Omega )} \left\| {\partial _{x_2 } u_{3m} }
\right\|_{L^2(\Omega )} \;+ } \right. \\
 &\qquad\qquad\qquad\qquad\qquad\qquad + \left. {\left\| {\partial _t u_{3m} }
\right\|_{L^4(\Omega )} \left\| {\partial _{x_3 } u_{3m} }
\right\|_{L^2(\Omega )} } \right) \\
 &\le \left( {\sum\limits_{i=1}^3 {\left\| {\partial _t u_{im} }
\right\|_{L^4(\Omega )}^2 } } \right)^{1/2}\left(
{\sum\limits_{j=1}^3 {\left\| {\partial _t u_{jm} }
\right\|_{L^4(\Omega )}^2 } } \right)^{1/2}\left(
{\sum\limits_{i,j=1}^3 {\left\| {\partial _{x_i } u_{jm} }
\right\|_{L^2(\Omega )}^2 } } \right)^{1/2} \\
\end{split}
\end{equation*}
where
\begin{equation*}
\begin{split}
 &\sum\limits_{i=1}^3 {\left\| {\partial _t u_{im} } \right\|_{L^4(\Omega )}^2 }
\le 2\sum\limits_{i=1}^3 {\left( {\left\| {\partial _t u_{im} }
\right\|_{L^2(\Omega )}^{1/2} \left\| {\nabla \partial _t u_{im} }
\right\|_{L^2(\Omega )}^{3/2} } \right)} \\
 &\quad \le 2\left( {\sum\limits_{i=1}^3 {\left\| {\partial _t u_{im} }
\right\|_{L^2(\Omega )}^2 } } \right)^{1/4}\left(
{\sum\limits_{i=1}^3 {\left\| {\nabla \partial _t u_{im} }
\right\|_{L^2(\Omega )}^2 } }
\right)^{3/4} \\
 \end{split}
\end{equation*}
so that
\begin{equation*}
\begin{split}
 &\partial _t \left( {\sum\limits_{i=1}^3 {\left\| {\partial _t u_{im} }
\right\|_{L^2(\Omega )}^2 } } \right)+2\left( {\sum\limits_{i=1}^3
{\left\|
{\nabla \partial _t u_{im} } \right\|_{L^2(\Omega )}^2 } } \right) \\
 &\quad \le 2^2\left( {\sum\limits_{i=1}^3 {\left\| {\partial _t u_{im} }
\right\|_{L^2(\Omega )}^2 } } \right)^{1/4}\left(
{\sum\limits_{i=1}^3 {\left\| {\nabla \partial _t u_{im} }
\right\|_{L^2(\Omega )}^2 } } \right)^{3/4}\left(
{\sum\limits_{i=1}^3 {\left\| {\nabla u_{im} }
\right\|_{L^2(\Omega )}^2 } } \right)^{1/2} \\
 &\quad \le 3^3 \left( {\sum\limits_{i=1}^3 {\left\| {\partial _t u_{im} }
\right\|_{L^2(\Omega )}^2 } } \right)\left( {\sum\limits_{i=1}^3
{\left\| {\nabla u_{im} } \right\|_{L^2(\Omega )}^2 } }
\right)^2 + \, \left( {\sum\limits_{i=1}^3 {\left\| {\nabla \partial _t
u_{im} } \right\|_{L^2(\Omega
)}^2 } } \right) \\
 \end{split}
\end{equation*}
it follows that
\[
\partial _t \left( {\sum\limits_{i=1}^3 {\left\| {\partial _t u_{im} }
\right\|_{L^2(\Omega )}^2 } } \right) +  \left( {\sum\limits_{i=1}^3
{\left\| {\nabla \partial _t u_{im} } \right\|_{L^2(\Omega )}^2 } }
\right)\le \phi_m (t)\left( {\sum\limits_{i=1}^3 {\left\| {\partial
_t u_{im} } \right\|_{L^2(\Omega )}^2 } } \right)
\]
where
\[
\phi_m (t)= 1\,+\, 3^3 \left( {\sum\limits_{i=1}^3 {\left\| {\nabla
u_{im} } \right\|_{L^2(\Omega )}^2 } } \right)^2
\]
\\

Introducing a stream function: $\psi =(\psi _2 ,\psi _2 ,\psi _3 )$,
\[
\mbox{curl}\psi =(\partial _{x_2 } \psi _3 -\partial _{x_3 } \psi _2
,\;\,\;\partial _{x_3 } \psi _1 -\partial _{x_1 } \psi _3
,\;\,\;\partial _{x_1 } \psi _2 -\partial _{x_2 } \psi _1 )
\]
According to $\omega =\mbox{curl}u$, $u=\mbox{curl}\psi $ and
$\mbox{div}\psi =0$, we have
\[
\mbox{curlcurl}\psi =-\Delta \psi =\omega , \quad -\Delta
\mbox{curl}\psi =\mbox{curl}\omega
\]
That is, $-\Delta u=\mbox{curl}\omega $. Then $(-\Delta
u,\;\,u)=(\mbox{curl}\omega ,\;\,u)$, where
\begin{equation*}
\begin{split}
 &(-\Delta u,\;\,u)=\sum\limits_{i=1}^3 {(-\Delta u_i ,\;\,u_i )}
=\sum\limits_{i=1}^3 {(\nabla u_i ,\;\,\nabla u_i )}
=\sum\limits_{i=1}^3
{\left\| {\nabla u_i } \right\|_{L^2(\Omega )}^2 } \\
 \end{split}
\end{equation*}
\begin{equation*}
\begin{split}
 &(\mbox{curl}\omega ,\;\,u)=(\partial _{x_2 } \omega _3 -\partial _{x_3 }
\omega _2 ,\;\;u_1 )+(\partial _{x_3 } \omega _1 -\partial _{x_1 }
\omega _3 ,\;\;u_2 )+(\partial _{x_1 } \omega _2 -\partial _{x_2 }
\omega _1 ,\;\;u_3
) \\
 &\quad \quad \quad =-(\omega _3 ,\;\partial _{x_2 } u_1 )+(\omega _2
,\;\partial _{x_3 } u_1 )-(\omega _1 ,\;\partial _{x_3 } u_2
)+(\omega _3 ,\;\partial _{x_1 } u_2 )-(\omega _2 ,\;\partial _{x_1
} u_3 )+(\omega _1
,\;\partial _{x_2 } u_3 ) \\
 &\quad \quad \quad =(\omega _1 ,\;\;\partial _{x_2 } u_3 -\partial _{x_3 }
u_2 )+(\omega _2 ,\;\;\partial _{x_3 } u_1 -\partial _{x_1 } u_3
)+(\omega
_3 ,\;\;\partial _{x_1 } u_2 -\partial _{x_2 } u_1 ) \\
 &\quad \quad \quad =(\omega ,\;\,\mbox{curl}u)=(\omega ,\omega
)=\sum\limits_{i=1}^3 {\left\| {\omega _i } \right\|_{L^2(\Omega )}^2 } \\
 \end{split}
\end{equation*}
Hence,
\begin{equation*}
\left( {\sum\limits_{i=1}^3 {\left\| {\nabla u_i }
\right\|_{L^2(\Omega )}^2 } } \right)^{1/2}=\left(
{\sum\limits_{i=1}^3 {\left\| {\omega _i } \right\|_{L^2(\Omega )}^2
} } \right)^{1/2}
\end{equation*}
it follows that
\[
\phi_m (t)= 1\,+\, 3^3 \left( {\sum\limits_{i=1}^3 {\left\| {\omega _{im} }
\right\|_{L^2(\Omega )}^2 } } \right)^2<+\infty
\]
\\

By the Gronwall inequality,
\[
\frac{d}{dt}\left\{ {\left( {\sum\limits_{i=1}^3 {\left\| {\partial
_t u_{im} } \right\|_{L^2(\Omega )}^2 } } \right)\;\exp \left(
{-\int_0^t {\phi_m (s)ds} } \right)} \right\}\le 0
\]
whence
\[
\mathop {\sup }\limits_{t\in (0,T)} \left( {\sum\limits_{i=1}^3
{\left\| {\partial _t u_{im} (t)} \right\|_{L^2(\Omega )}^2 } }
\right)\le \left( {\sum\limits_{i=1}^3 {\left\| {\partial _t u_{im}
(0)} \right\|_{L^2(\Omega )}^2 } } \right)\;\exp \left( {\int_0^T
{\phi_m (s)ds} } \right)
\]
Therefore
\[
\partial _t u_{im} \in L^\infty (0,T;\;H) \cap L^2 (0,T;\;V), \qquad \qquad  i=1,2,3
\]
\\

Similar to the Theorem 3.8 in Chapter 3 of [4], we obtain
\[
u_i \in L^\infty (0,T;\;H^2(\Omega )),  \qquad \qquad i=1,2,3
\]
\\

Noting that $(-\Delta u,\;v)=(-\partial _t
u-(u\cdot \nabla )u,\;\,v)$. Since $\partial _t u$ and $(u\cdot \nabla
)u$ are of some degree of continuity, then $u$ can reach a higher
degree of continuity, based on the smoothing effect of inverse
elliptic operator $\Delta ^{-1}$. By repeated application of this
process one can prove that the solution $u$ is in $C^\infty (\Omega
\times (0,T))$.
\\
\\
\\

\end{document}